\documentclass[11pt]{article}
\usepackage[a4paper,top=3cm,bottom=2cm,left=3cm,right=3cm,marginparwidth=1.75cm]{geometry}
\usepackage{graphicx}%
\usepackage{multirow}%
\usepackage{amsmath,amssymb,latexsym,amstext,amsthm}
\usepackage{amsfonts}%
\usepackage{mathrsfs}%
\usepackage{hyperref}

\theoremstyle{thmstyleone}%
\newtheorem{theorem}{Theorem}[section]        
\newtheorem{lemma}{Lemma}[section]
\newtheorem{corollary}{Corollary}[section]
\newtheorem{definition}{Definition}[section]

\theoremstyle{thmstyletwo}%

\newcommand{\Pochh}[2]{\left( #1 \right)_{#2}}
\newcommand{\supp}{\text{supp }}
\providecommand{\keywords}[1]
{
  \small	
  \textbf{\textit{Keywords:}} #1
}

\providecommand{\subjclass}[1]
{\small	
  \textbf{\textit{2010 Mathematical Subject Classification:}} #1
}

\title{On two families of discrete multiple orthogonal polynomials on a
star-like set\footnote{This is a preprint accepted for publication in Lobachevskii Journal of Mathematics, 2023, No. 12, \url{https://link.springer.com/journal/12202} \copyright Pleiades Publishing, Ltd., 2023.}}

\author{Jorge Arves\'u Carballo\\
\small Department of Mathematics, Universidad Carlos III de Madrid, Avda. de la
Universidad, 30 \\
\small Legan\'es, 28911 Madrid, Spain\\
Alejandro J. Quintero Roba\\
\small Department of Mathematics, Baylor University, 1410 S
4th Street. Waco, 76706, Texas, USA}

\date{December 30, 2022}

\begin{document}
\maketitle

\begin{abstract} We study two families of type II discrete multiple
orthogonal polynomials on an $r$-legged star-like set with respect
to $r$ weight functions of Charlier (Poisson distributions) and
Meixner (negative binomial distributions), respectively. We focus
our attention on the structural properties such as the raising
operators, the Rodrigues-type formulas, and the explicit expressions
of the polynomial families as well as the coefficients of the
nearest neighbor recurrence relation. A limit relation between these
families is given.
\end{abstract}

\noindent\subjclass{33C45, 33C47, 41A28, 42C05}\\
\noindent\keywords{Multiple orthogonal polynomials, Charlier--Angelesco
polynomials, Meix\-ner--An\-ge\-les\-co polynomials, Rodrigues-type formula,
Recurrence relation}

\section{Introduction}
Multiple orthogonal polynomials are an  extension of orthogonal
polynomials where the orthogonality conditions are shared with
respect to a system of measures
\cite{bib:aptekarev,bib:kalyagin,bib:kalyagin-ronveaux,bib:nikishin-sorokin}.
They appeared in connection with the simultaneous rational
approximation of a system of analytic functions
\cite{bib:Gonchar,bib:kn_Nikishin,bib:Sorokin3,bib:sorokin1}. This
notion of simultaneous approximation traces back to Hermite's proof
\cite{Hermite_exp} of the transcendence of the number $e$. There are
two types of multiple orthogonal polynomials (type I and type II).
See survey \cite{bib:aptekarev} and Chapter 4 of the monograph
\cite{bib:nikishin-sorokin} for a detailed information about
simultaneous Pad\'e approximation and the two types of multiple
orthogonal polynomials.

Consider a multi-index $\vec{n}\in\mathbb{Z}_{+}^{r}$ of length $\lvert\vec{n}\rvert=n_0+\cdots+n_{r-1}$ and the vector measure $\vec{\mu}=(\mu_{0},\ldots,\mu_{r-1})$, where the components are positive measures with finite moments. For the multi-index $\vec{n}$, a monic polynomial $P_{\vec{n}}$ of degree $\lvert\vec{n}\rvert$ is said to be a type II multiple orthogonal polynomial if the following orthogonality relations hold
\begin{equation}
\int P_{\vec{n}}(x)x^k\ d\mu_{\ell}(x) =  0,\quad k=0,\ldots,n_{\ell}-1,\quad \ell=0,\ldots,r-1.
\label{eq:def_typeII}
\end{equation}
These relations lead to a linear system of $\lvert\vec{n}\rvert$ equations for $\lvert\vec{n}\rvert$ unknown coefficients of $P_{\vec{n}}$. If there exists a unique solution $P_{\vec{n}}$, we call $\vec{n}$ a normal index.
The existence and uniqueness of type II and type I multiple orthogonal polynomials are equivalent \cite{bib:nikishin-sorokin} since the matrix of coefficients derived from \eqref{eq:def_typeII} is the transpose of the corresponding matrix for the type I multiple orthogonal polynomials. In this paper we will only deal with type II multiple orthogonal polynomials $P_{\vec{n}}$ with complex measures. If all multi-indices $\vec{n}$ are normal we say that the system of measures $\{\mu_{0},\ldots,\mu_{r-1}\}$ is a perfect system \cite{bib:Mahler,bib:nikishin-sorokin}. Two important examples of perfect systems of measures are the AT-system \cite{bib:nikishin-sorokin} and Angelesco system \cite{bib:angelesco,bib:nikishin}. In an AT-system all the components of $\vec{\mu}$ have the same support (see \cite{bib:arvesu-coussement-vanassche,bib:vanassche-coussement} for examples involving classical measures), while in an Angelesco system, each component $\mu_{\ell}$ has support $\Omega_{\ell}$, where $\Omega_0,\ldots,\Omega_{r-1}$ are disjoint and are allowed to touch, i.e. $\Omega_i^{\circ} \cap \Omega_j^{\circ} = \emptyset$,  for  $i\neq j$ (see \cite{bib:lee,bib:sorokin2} for some examples that also involve classical measures).

The type II multiple orthogonal polynomial $P_{\vec{n}}(z)$ from \eqref{eq:def_typeII} is the common denominator of the simultaneous rational approximants $\frac{Q_{\vec{n},\ell}(z)}{P_{\vec{n}}(z)}$ to the Cauchy transforms
\begin{equation*}
\displaystyle \hat{\mu}_{\ell}(z)=\int\frac{d\mu_{\ell}(x)}{z-x},
\quad z\notin\mbox{supp}\,\mu_{\ell},\quad \ell=0,\dots,r-1,
\end{equation*}
of the vector components of $\vec{\mu}=(\mu_{0},\ldots,\mu_{r-1})$, that is, the following simultaneous rational approximation with prescribed order near infinity holds
\begin{equation*}
\begin{array}{c}
\displaystyle P_{\vec{n}}(z)\hat{\mu}_{\ell}(z)-Q_{\vec{n},\ell}(z)=\frac{\zeta_{\ell}}{z^{n_{\ell}+1}}+\cdots=\mathcal{O}(z^{-n_{\ell}-1}),\quad
\ell=0,\dots,r-1.
\end{array}
\label{f1}
\end{equation*}

If the measures in~\eqref{eq:def_typeII} are discrete
\begin{gather}
\mu_{\ell}=\sum\limits_{k=0}^{N_{\ell}}\omega_{\ell,k}\delta_{x_{\ell,k}},
\quad
\omega_{\ell,k}>0,
\qquad
N_{\ell}\in \mathbb{N\cup}\{+\infty \},
\qquad
\ell=0,\ldots,r-1,
\notag
\end{gather}
where $\delta_{x_{\ell,k}}$ denotes the Dirac delta function and $x_{\ell_{1},k}\neq x_{\ell_{2},k}$, $k=0,\ldots,N_{\ell}$,
whenever $\ell_{1}\neq \ell_{2}$, we say that $P_{\vec{n}}(x)$ is a discrete multiple orthogonal polynomial \cite{bib:arvesu-coussement-vanassche,bib:arvesu-ramirez1,bib:arvesu-ramirez2,bib:arvesu-ramirez3}. Two interesting examples are the multiple Charlier polynomials and multiple Meixner polynomials of the first kind (see \cite{bib:miki-tsujimoto-vinet-zhedanov1,bib:miki-tsujimoto-vinet-zhedanov2,bib:ndayiragije1-vanassche} for a mathematical approach to a physical phenomenon involving these polynomials).

The type II multiple Charlier polynomials $C_{\vec{n}}^{\vec{a}}(x)$
are orthogonal with respect to an AT-system of discrete measures
supported on the linear lattice $ x_k = k$,
$k\in\mathbb{N}\cup\{0\}$, with weight functions $\omega_\ell(x) =
{a_{\ell}^{x}}/{\Gamma(x+1)},$ where $a_\ell > 0$ are pairwise
different. These polynomials have the Rodrigues-type formula
\begin{equation}\label{eq:ch_mclassic}
C_{\vec{n}}^{\vec{a}} (x) = \prod\limits_{j=0}^{r-1} (-a_j)^{n_j} \Gamma(x+1) \prod\limits_{j=0}^{r-1} \left( a_j^{-x} \bigtriangledown ^{n_j} a_j^x \right) \dfrac{1}{\Gamma(x+1)},
\end{equation}
which is used to obtain an explicit expression of the polynomials and deduce a nearest neighbor recurrence relation \cite{bib:arvesu-coussement-vanassche,bib:vanassche}. Here $\bigtriangledown$ is the standard backward difference operator $\bigtriangledown f(x)=f(x)-f(x-1)$. In the sequel, we will use a distinct backward difference operator $\nabla$ acting on composite functions as follows
\begin{equation}
\nabla f(s(z)) = f(s(z)) - f(s(z)-1),\quad
\nabla f(s(z)+1)=\Delta f(s(z)),\label{relation_delta_nabla}
\end{equation} with $s(z)\pm 1$ also in the domain of $f$.
Moreover, the following relations will be used when convenient:
\begin{align}
\nabla f(s(z))g(s(z)) &= f(s(z)) \nabla g(s(z)) + g(s(z)-1) \nabla f(s(z)),\notag\\
\nabla^n f(s(z)) &= \displaystyle{\sum\limits_{k=0}^{n} {n\choose{k}}}\left(-1\right)^k f(s(z)-k),\label{nth_nabla}\\
\nabla \left( s(z) \right)_{j+1} &= (j+1) \left( s(z) \right)_j,\quad j \in \mathbb{N},\notag
\end{align}
where $\Pochh{s(z)}{j+1}=s(z)(s(z)+1)\ldots(s(z)+j)$ and $\left( s(z) \right)_0=1$, denotes the Pochhammer symbol.

The type II multiple Meixner polynomials of the first kind $M_{\vec{n}}^{\beta, \vec{c}}(x)$ are orthogonal with respect to an AT-system of discrete measures supported on the linear lattice $x_k = k$, with weight functions $$\omega_\ell(x) = \dfrac{\Gamma(x+\beta)}{\Gamma(\beta)\Gamma(x+1)}c_{\ell}^{x},$$ where $\beta \neq 0,-1,-2,\ldots$ and $0< c_\ell < 1$ pairwise different. These polynomials have the Rodrigues-type formula
\begin{equation}\label{eq:m1_mclassic}
M_{\vec{n}}^{\beta,\vec{c}} (x) = \prod\limits_{j=0}^{r-1} \left( \dfrac{c_j}{c_j-1} \right)^{n_j} \dfrac{\Gamma(x+1) }{\Gamma(x+\beta)}\prod\limits_{j=0}^{r-1} \left( c_j^{-x} \bigtriangledown ^{n_j} c_j^x \right) \dfrac{\Gamma(x+ \lvert\vec{n}\rvert+\beta)}{\Gamma(x+1)},
\end{equation}
and satisfy a nearest neighbor recurrence relations \cite{bib:arvesu-coussement-vanassche,bib:vanassche}.

Recently in \cite{bib:leurs-vanassche2,bib:leurs-vanassche1} some
examples of continuous multiple  orthogonal polynomials with respect
to Angelesco systems of measures on a star-like set were considered:
The type I Jacobi--Angelesco and the types I and II
Laguerre--Angelesco polynomials, respectively. In particular, for
the Laguerre--Angelesco polynomials the measures $\mu_\ell$,
$\ell=0,\ldots,r-1$, are supported on an $r$-star formed with
intervals $\Omega_\ell = [0,e^{2\pi i\ell /r}\infty)$,
$\ell=0,\ldots,r-1$, and weight functions
$
w(z) = \lvert z\rvert^\beta e^{-z^r}, \ \ \beta>-1.
$
The above $r$-star was already considered in
\cite{bib:aptekarev-kalyagin-saff,bib:delvaux-lopez},  however the
authors dealt with Nikishin systems of measures, which are also
perfect systems. The type II Laguerre--Angelesco orthogonal
polynomials is the monic polynomial $L_{\vec{n}}^{\beta}(z)$, with
degree exactly $\lvert\vec{n}\rvert$ in $z^{r}$, defined by the
orthogonality conditions
\begin{equation*}
\int\limits_{\Omega_\ell} z^k L_{\vec{n}}^{\beta} (z) \lvert z\rvert^{\beta} e^{-z^r} \text{d}z = 0,\quad 0\leq k \leq n_\ell-1,\quad \ell= 0, \ldots, r-1.
\end{equation*}
For the diagonal multi-index $\vec{n}= (n, n,\ldots, n)$, it was proven in \cite{bib:leurs-vanassche1} that these polynomials verify a Rodrigues-type formula and a nearest neighbor recurrence relation, among other algebraic properties.

Motivated by \cite{bib:leurs-vanassche1} as well as
\cite{bib:arvesu-coussement-vanassche}, where among some discrete
multiple orthogonal polynomials the multiple Charlier and Meixner
polynomials of the first kind were studied, we set as a goal of this
contribution the investigation of Charlier--Angelesco and
Meixner--Angelesco polynomials with respect to $r$ complex measures
supported on an $r$-star defined by the intervals $[0,e^{2\pi i
k/r}\infty)$, $k=0,\ldots, r-1$, with finite moments.

The paper is structured as follows. In Section \ref{sec:2}  we
introduce the notion of discrete multiple orthogonal polynomials on
an $r$-star and address the location of zeros as well as the
recurrence relation. In Section \ref{sec:ch_m1_measures} we define
the systems of measures that will be used in Sections
\ref{sec:ch2_examples} and \ref{sec:meix2_examples}. These sections
contain the main results. We investigate two families of type II
discrete multiple orthogonal polynomials on an $r$-star, namely the
Charlier--Angelesco and Meixner--Angelesco polynomials of the first
kind. For these polynomials, we obtain the raising operators,
Rodrigues-type formulas, and their explicit expressions. With these
expressions we compute the recurrence coefficients of the nearest
neighbor recurrence relations. We end with some concluding remarks
in Section \ref{conclusions} and give a limit relation involving the
Charlier--Angelesco and Meixner--Angelesco polynomials.

\section{Type II discrete multiple orthogonal polynomials on an $r$-star}\label{sec:2}\

Let $\omega = e^{2\pi i/r}$. For $r\in \mathbb{N}$ and
$j=0,1,\ldots, r-1$, $\omega^j$ are the $r$th roots of $1$. Define
the intervals $\Omega_j$ as counterclockwise rotations of the
positive real axis,
$
\Omega_j =  \left[ 0, \omega^j\infty \right) = \omega^j \cdot \mathbb{R_+}\cup \{ 0 \}.
$
With these intervals we form the $r$-star on the complex plane
\begin{equation}\label{r-star}
\bigcup\limits_{\ell=0}^{r-1}\Omega_\ell.
\end{equation}

On each interval $\Omega_j$, $j=0,\ldots,r-1$, we will consider a continuous function $\rho_j$ (weight function) and a discrete measure on the mass points $z_{j,0},\ldots,z_{j,N}\in \Omega_{j}$,
\begin{equation} \label{eq:particular_measures}
\mu_j
=\sum_{k=0}^{N} \rho_{j}\delta_{z_{j,k}},\quad N \in \mathbb{N} \cup \{ \infty \},
\end{equation}
with finite moments $m_{\ell}^{(j)}$ $(\ell = 0,1,2, \ldots)$ and full column rank matrix
\begin{eqnarray}
\label{matrix}
\mathcal{M} & = & \left(
\begin{array}{c}
\mathcal{M}_0(n_1)\\
\vdots\\
\mathcal{M}_{r-1}(n_r)
\end{array} \right),\quad \mathcal{M}_j(n_{j+1}) = \left(
\begin{array}{cccc}
m_0^{(j)} & m_1^{(j)} & \ldots & m_{\lvert\vec{n}\rvert}^{(j)}\\
\vdots & \vdots & & \vdots\\
m_{n_j-1}^{(j)} & m_{n_j}^{(j)} & \ldots &
m_{\lvert\vec{n}\rvert+n_j-1}^{(j)}
\end{array} \right),
\end{eqnarray}
where $\mathcal{M}_j(n_{j+1})$ is the
$n_{j}\times(\lvert{\vec{n}}\rvert+1)$  matrix of moments of the
measure $\mu_j$, for $\vec{n}\in\mathbb{N}^{r}$.

Next, we define sequences of discrete multiple orthogonal
polynomials in  powers of $z^r$ on an $r$-star and set $\supp
{\mu_j} = \Omega_j$, with $\Omega_i^{\circ} \cap \Omega_j^{\circ} =
\emptyset$,  for  $i\neq j$.

\begin{definition} \label{def:typeII}
The type II discrete multiple orthogonal polynomial $P_{\vec{n}}$
for  the multi-index $\vec{n}$ on the $r$-star \eqref{r-star} is the
polynomials of degree $\leq r\lvert\vec{n}\rvert$ defined by the
orthogonality relations
\begin{align}
\int\limits_{\Omega_0} (z^r)^k P_{\vec{n}} (z) \text{d} \mu_0 &= 0, \ \ \ k = 0, \ldots, n_0-1, \nonumber\\
&\vdots \label{eq:type2_original}\\
\int\limits_{\Omega_{r-1}} (z^r)^k P_{\vec{n}} (z) \text{d} \mu_{r-1} &= 0, \ \ \ k = 0, \ldots, n_{r-1}-1. \nonumber
\end{align}
\end{definition}
These orthogonality conditions are equivalently written as
\begin{align}\label{eq:discrete_orth_conditions}
\displaystyle\sum_{k=0}^{\infty} P_{\vec{n}}(z_{\ell,k})
\Pochh{-z_{\ell,k}^r}{j}  \rho_{\ell}(z_{\ell,k}) &= 0,\quad 0\leq
\ell  \leq r-1,\quad 0\leq j\leq n_{\ell}-1,
\end{align}
where $\Pochh{-z^r}{n_{\ell}-1}=(-z^r)(1-z^r)\cdots(n_{\ell}-2-z^r)$
is the Pochhammer symbol, i.e., a polynomial of degree
$r(n_{\ell}-1)$. Here, instead of the canonical basis in the vector
subspace of polynomials we have used the basis of Pochhammer
symbols.

Observe that the orthogonality conditions \eqref{eq:type2_original} define a linear systems of $\lvert\vec{n}\rvert$ equations for $\lvert\vec{n}\rvert+1$ known coefficients of $P_{\vec{n}}$, with $\mathcal{M}$ as the matrix of this linear system. We focus on the solution $P_{\vec{n}}$ that is unique up to a multiplicative factor and also on polynomial of degree exactly $\lvert\vec{n}\rvert$ (the monic polynomial exists and will be unique). In this situation, $\vec{n}$ is a normal index. In what follows, we will only deal with normal indices $\vec{n}\in\mathbb{N}^{r}$, that is, when the matrix
$\mathcal{N}$ given from $\mathcal{M}$ by deleting the last column has rank $\lvert\vec{n}\rvert$. Thus, the corresponding unique monic multiple orthogonal polynomial solution of \eqref{eq:type2_original} is
\begin{equation}
\displaystyle P_{\vec{n}}(z) = \frac{1}{\det \mathcal{N}} \det\left(
\begin{array}{cccc}
m_0^{(0)} & m_1^{(0)} & \ldots & m_{\lvert\vec{n}\rvert}^{(0)}\\
\vdots & \vdots & & \vdots\\
m_{n_1-1}^{(0)} & m_{n_1}^{(0)} & \ldots &
m_{\lvert\vec{n}\rvert+n_1-1}^{(0)}\\
\vdots & \vdots & & \vdots\\
m_0^{(r-1)} & m_1^{(r-1)} & \ldots & m_{\lvert\vec{n}\rvert}^{(r-1)}\\
\vdots & \vdots & & \vdots\\
m_{n_r-1}^{(r-1)} & m_{n_r}^{(r-1)} & \ldots &
m_{\lvert\vec{n}\rvert+n_r-1}^{(r-1)}\\ 1 & z^{r} & \ldots &
z^{r\lvert\vec{n}\rvert}
\end{array} \right),\label{Pn_matrix_eq}
\end{equation}
where
\begin{equation}\label{eq:particular_moments}
m_{\ell}^{(j)} = \int\limits_{\Omega_j} (z^r)^{\ell}\text{d}\mu_j = \sum_{k=0}^{N} \left( z_{j,k}^r \right)^{\ell} \rho_j (z_{j,k}^r),\quad j=0,\ldots,r-1,\quad \ell = 0,1, \ldots
\end{equation}

Now, we address a couple of algebraic results of type II  Angelesco
polynomials on an $r$-star; namely, the zero location and recurrence
relation theorems, respectively.

\begin{definition}\label{def:w-symmetric}
The system of weight functions $\{ \rho_{0}(z),\rho_1(z)\ldots,
\rho_{r-1}(z)\}$ is said to be $\omega$-symmetric if for every
$j=0,\ldots,r-1$ and $k\in\mathbb{N}$
$
\rho_j(\omega^{k}z)=\rho_j(z).
$
\end{definition}

Notice that this $\omega$-symmetry is inherited by the system of
measure \eqref{eq:particular_measures}. Moreover, in
\eqref{eq:type2_original}, the orthogonality relations
$$
\displaystyle\int_{0}^{\omega^{j}\infty} (z^r)^k P_{\vec{n}}(z) \text{d} \mu_j = 0, \quad k = 0, \ldots, n_j-1, \quad j=1,\ldots,r-1,
$$
can be interpreted as rotated copies of the relation
$$
\displaystyle\int_{0}^{\infty} (z^r)^k P_{\vec{n}} (z) \text{d} \mu_0 = 0, \quad k = 0, \ldots, n_0-1,
$$
because the integrand expression has an $\omega$-symmetry, since $P_{\vec{n}}$ is a polynomial in $z^r$.

\begin{theorem}\label{location_zeros} Let $\vec{n}=(n,n,\ldots,n)$ and $\rho_{j}(z)$ be the $\omega$-symmetric weight functions supported on $\Omega_j$ $(j=0,\ldots,r-1)$ of the $r$-star \eqref{r-star}. Then, all the zeros of $P_{\vec{n}}$ in \eqref{Pn_matrix_eq} are simple and lie in the rays of the $r$-star.
\end{theorem}

We will see that these zeros are rotated copies of the real zeros in the interval $(0,\infty)$.

\begin{proof} Assume $n>m\in\mathbb{N}$ and $x_1,\ldots, x_m$
are the zeros of odd multiplicity of $P_{\vec{n}}$ that belong to the interval $(0,\infty)$. Consider a polynomial $Q(z)$ of degree $rm$ with exactly $m$ sign changes on $(0,\infty)$ at the points $x_1,\ldots, x_m$ and no other zeros. Set $Q(z)=(z^r-x_{1}^{r})\cdots(z^r-x_{m}^{r})$ then,
the product $Q(z)P_{\vec{n}}(z)$ does not change sign on $(0,\infty)$. Hence,
$$
\int_{0}^{\infty} Q(z)P_{\vec{n}} (z) \text{d} \mu_0 \not= 0.
$$
However, by the first orthogonality relation in
\eqref{eq:type2_original}  this integral is zero. This contradiction
implies that $m\geq n$. Moreover, $P_{\vec{n}}$ is a polynomial of
degree exactly $|\vec{n}|$ in $z^r$ and, for every non-negative
integer $\ell$, it verifies the $\omega$-symmetry property
$P_{\vec{n}}(\omega^{\ell}z)=P_{\vec{n}}(z)$, that is, it has the
same number of zeros in each ray, which implies that $m=n$.
Therefore, the zeros on the rays of the $r$-star are rotated copies
of the zeros $x_1, \ldots, x_{n}$ lying in the interval
$(0,\infty)$.
\end{proof}

Now we investigate some preliminary results used in the next
Sections \ref{sec:ch2_examples} and \ref{sec:meix2_examples}.

\begin{lemma}\label{lemma:Pn_basis}
Let $\mathbb{V}$ be the vector subspace of polynomials $Q(z)$ in the variable $z^r$ of degree at most $\lvert\vec{n}\rvert-1$ defined by the following conditions
\begin{equation}\label{eq:Pn_lemma_basis}
\int\limits_{\Omega_\ell} Q(z) \Pochh{-z^r}{j} \text{d} \mu_{\ell} = 0, \quad 0\leq j\leq n_{\ell}-2,\quad 0\leq \ell\leq r-1.
\end{equation}
Then, the system $\left\lbrace P_{\vec{n}-\vec{e}_{\ell}}\right\rbrace_{\ell=0}^{r-1}$ of multiple orthogonal polynomials of degree $\lvert\vec{n}\rvert-1$ with respect to the system $\left\lbrace \mu_0, \ldots, \mu_{r-1}\right\rbrace$ is a basis of $\mathbb{V}$.
\end{lemma}

Here, the multi-index $\vec{e}_{\ell}$ denotes  the standard
$r$-dimensional unit vector with the $\ell$th entry equals $1$ and
$0$ otherwise.

\begin{proof} From the orthogonality relations \eqref{eq:type2_original} we have
that $P_{\vec{n}-\vec{e}_{\ell}}$, for every $\ell=0,\ldots,r-1$,
satisfies \eqref{eq:Pn_lemma_basis}. Hence, it belongs to
$\mathbb{V}$. Now, we prove that these polynomials form a basis of
$\mathbb{V}$, which implies that $\dim \mathbb{V} = r$. Indeed, we
implement a proof by contradiction. Suppose $\left\lbrace
P_{\vec{n}-\vec{e}_{\ell}}\right\rbrace_{\ell=0}^{r-1}$ is a
linearly dependent system. Hence, there exists some numbers
$\lambda_\ell$ $(\ell= 0,\ldots,r-1)$ such that
\begin{equation}\label{eq:ch_ld_condition}
\displaystyle{\sum_{\ell=0}^{r-1}} \lambda_\ell P_{\vec{n}-\vec{e}_{\ell}} (z) = 0,\quad
\mbox{where}\quad\sum_{\ell=0}^{r-1}\lvert\lambda_\ell\rvert>0.
\end{equation}

Now, for a fixed $j$, where $0\leq j  \leq r-1$, multiply  the
previous equation by $\Pochh{-z^r}{n_j-1} \mu_j$ and then integrate
over the $j$th ray
\begin{align}\label{lambda_integ}
\displaystyle{\sum_{\ell=0}^{r-1}} \lambda_\ell
\int_{0}^{\omega^{j}\infty}P_{\vec{n}-\vec{e}_{\ell}}(z)
\Pochh{-z^r}{n_j-1} \text{d}\mu_j &= 0.
\end{align}

From the orthogonality relations \eqref{eq:type2_original} we have that
$$
\displaystyle\int_{0}^{\omega^{j}\infty}P_{\vec{n}-\vec{e}_{\ell}}(z) \Pochh{-z^r}{n_j-1} \text{d}\mu_j
= \left\{\begin{array}{l} =0,\quad\mbox{for}\quad \ell\not= j,\\\not=0,\quad\mbox{otherwise}. \end{array}\right.
$$
Therefore, the expression \eqref{lambda_integ} becomes
$$
\lambda_j \int_{0}^{\omega^{j}\infty}P_{\vec{n}-\vec{e}_{j}}(z) \Pochh{-z^r}{n_j-1} \text{d}\mu_j = 0,
$$
which implies that $\lambda_j=0$. Hence, repeating the same argument
for all other $j=0,\ldots,r-1$,  the coefficients $\lambda_j$
$(0\leq j\leq r-1)$ in \eqref{eq:ch_ld_condition} vanish. This
contradiction implies that $\left\lbrace
P_{\vec{n}-\vec{e}_{\ell}}\right\rbrace_{\ell=0}^{r-1}$ is a
linearly independent system in $\mathbb{V}$. Since any polynomial
from $\mathbb{V}$ has at most $\lvert\vec{n}\rvert$ coefficients and
$\mathbb{V}$ is defined by $(\lvert\vec{n}\rvert-r)$ equations
\eqref{eq:Pn_lemma_basis}, then $\dim \mathbb{V} \leq r$. But,
$\left\lbrace
P_{\vec{n}-\vec{e}_{\ell}}\right\rbrace_{\ell=0}^{r-1}\subset\mathbb{V}$
is a linearly independent system of vectors, so it spans
$\mathbb{V}$, that is, a basis of the vector subspace $\mathbb{V}$,
where $\dim \mathbb{V} = r$.
\end{proof}

The following notation for the coefficients of the monic polynomials $P_{\vec{n}}$ in $z^r$ is used:
\begin{equation}\label{monic_Pns}
P_{\vec{n}}(z) = \sum\limits_{m=0}^{\lvert\vec{n}\rvert} \alpha_{m}^{\vec{n}} z^{rm},
\quad
 P_{\vec{n}+\vec{e}_k}(z) = \sum\limits_{m=0}^{\lvert\vec{n}\rvert+1} \alpha_{m}^{\vec{n}+\vec{e}_k} z^{rm},\\
\end{equation}
where $\alpha_{\lvert\vec{n}\rvert}^{\vec{n}}
= \alpha_{\lvert\vec{n}\rvert+1}^{\vec{n}+\vec{e}_k} = 1$.

Now, we investigate the recurrence relation involved the type II discrete multiple orthogonal polynomials on the $r$-star (see \eqref{eq:recurrence_relation} below).

\begin{theorem}\label{theorem:recurrence_relation}
Suppose all the multi-indices $\vec{n}\in\mathbb{N}^{r}$ are normal
for  the system of measures $\mu_0, \ldots, \mu_{r-1}$. Then, the
corresponding type II discrete multiple orthogonal polynomials on
the $r$-star satisfy the recurrence relation
\begin{equation}\label{eq:recurrence_relation}
\left(z^r - b_{\vec{n},k} \right)P_{\vec{n}}(z)  - P_{\vec{n}+\vec{e}_k}(z) = \sum_{\ell=0}^{r-1} d_{\vec{n},\ell}  P_{\vec{n}-\vec{e}_\ell}(z),
\quad k=0,\ldots, r-1,
\end{equation}
 where
\begin{eqnarray}\label{eq:rr_linear_combination1}
b_{\vec{n},k}&=& \alpha_{\lvert\vec{n}\lvert-1}^{\vec{n}} -
\alpha_{\rvert\vec{n}\lvert}^{\vec{n}+\vec{e}_k}, \quad \label{eq:b_linear_combination1}\\
d_{\vec{n},\ell} = \dfrac{\displaystyle{\int_{\Omega_{\ell}}} z^r P_{\vec{n}}(z)
\Pochh{-z^r}{n_{\ell}-1}\text{d} \mu_\ell}{\displaystyle{\int_{\Omega_{\ell}}}
P_{\vec{n}-\vec{e}_{\ell}}(z)\Pochh{-z^r}{n_{\ell}-1}\text{d} \mu_\ell},
\quad \ell=0,\ldots,r-1. \notag
\end{eqnarray}
\end{theorem}

The proof differs slightly from that used for the  recurrence
relation of standard orthogonal polynomials (see, e.g., \cite{nu})
except for the use of Lemma \ref{lemma:Pn_basis}. Despite this fact
we prove it here.

\begin{proof} Since we are dealing in \eqref{eq:recurrence_relation}
with monic polynomials \eqref{monic_Pns}, we have that $z^r
P_{\vec{n}} - P_{\vec{n}+\vec{e}_k}$  is a polynomial of degree
$\lvert\vec{n}\rvert$ in $z^r$. Hence, taking $ b_{\vec{n},k}=
\alpha_{\lvert\vec{n}\rvert-1}^{\vec{n}} -
\alpha_{\lvert\vec{n}\rvert}^{\vec{n}+\vec{e}_k}, $ in the left-hand
side of the equation \eqref{eq:recurrence_relation}, the  resulting
polynomial $Q(z) = z^r P_{\vec{n}}(z) - P_{\vec{n}+\vec{e}_k}(z) -
b_{\vec{n},k}P_{\vec{n}}(z)$ has degree $\lvert\vec{n}\rvert-1$ in
$z^r$. Indeed, from the orthogonality relations
$$\displaystyle{\int\limits_{\Omega_\ell}} Q(z)\Pochh{-z^r}{j} d\mu_\ell = 0, \ \ \ j=0,\ldots, n_{\ell}-2, \ \ \ell=0,\ldots, r-1,
$$
we conclude that $Q(z)$ belongs to the subspace $\mathbb{V}$ defined
in   Lemma \ref{lemma:Pn_basis}. Therefore,
\begin{equation*}
    Q(z)=\left(z^r- b_{\vec{n},k}\right)P_{\vec{n}}(z) -
    P_{\vec{n}+\vec{e}_k}(z) = \sum\limits_{j=0}^{r-1} d_{\vec{n},j} P_{\vec{n}-\vec{e}_j}(z),
\end{equation*}
Multiplying \eqref{eq:recurrence_relation} by
$\Pochh{-z^r}{n_{\ell}-1}$  and integrating over the ray
${\Omega_\ell}$ with respect to the measure $\mu_{\ell}$,
$$
\int\limits_{\Omega_\ell}\left[\left(z^r - b_{\vec{n},k}
\right)P_{\vec{n}}(z)  -
P_{\vec{n}+\vec{e}_k}(z)\right]\Pochh{-z^r}{n_{\ell}-1}
\text{d}\mu_\ell =\sum_{j=0}^{r-1}
d_{\vec{n},j}\int\limits_{\Omega_\ell}
P_{\vec{n}-\vec{e}_j}(z)\Pochh{-z^r}{n_{\ell}-1}\text{d}\mu_\ell.
$$
From the orthogonality conditions \eqref{eq:type2_original} one gets
\begin{equation*}
\int\limits_{\Omega_\ell} z^r P_{\vec{n}}(z)\Pochh{-z^r}{n_{\ell}-1}\text{d}\mu_\ell = d_{\vec{n},\ell}\int\limits_{\Omega_{\ell}}    P_{\vec{n}-\vec{e}_{\ell}}(z)\Pochh{-z^r}{n_{\ell}-1}\text{d}\mu_\ell, \end{equation*}
equivalently,
$$
d_{\vec{n},\ell} = {\displaystyle{\int_{\Omega_\ell}} z^r
P_{\vec{n}}(z)\Pochh{-z^r}{n_{\ell}-1}
\text{d}\mu_\ell}/{\displaystyle{\int_{\Omega_\ell}}
P_{\vec{n}-\vec{e}_{\ell}}(z)\Pochh{-z^r}{n_{\ell}-1}\text{d}\mu_\ell}.
$$
\end{proof}

\section{On two systems of discrete measure on an $r$-star}\label{sec:ch_m1_measures}

Here we introduce two systems of measures that  will be used in
Sections \ref{sec:ch2_examples} and \ref{sec:meix2_examples},
respectively. There we will find the raising operators,
Rodrigues-type formulas and explicit expressions for the
corresponding families of Angelesco multiple orthogonal polynomials.
In addition, at the end of these sections, the coefficients for the
nearest neighbor recurrence relations in Theorem
\ref{theorem:recurrence_relation} will be calculated.

\subsection{Charlier--Angelesco system of measures}
Let $\vec{a}=(a_0,\ldots,a_{r-1})\in \mathbb{C}^r$ be a vector with non-zero complex entries. Moreover, we consider the mass points
\begin{equation} \label{eq:mass_points}
\{z_{j,k}\,\,:\,\,z_{j,k} = k^{1/r}\omega^j\}_{k=0}^\infty \subset \Omega_j,\quad j=0,\ldots,r-1,
\end{equation}
and the weight functions,
\begin{equation}\label{eq:ch_weight}
\rho_j (z) = \dfrac{a_j^{z^{r}}}{\Gamma(z^{r}+1)}=\dfrac{e^{z^{r}\ln a_j}}{\Gamma(z^{r}+1)}, \ \ \ j=0,\ldots, r-1,
\end{equation}
where $-\pi<\arg a_j\leq \pi$ (principal branch of the logarithmic function). The function $\rho_j (z)$ is an extension of the Poisson distribution $\rho(k)=a^{k}/k!$, $a>0$, on the non-negative integers $k=0,1,2\ldots$, where a complex replacement of the real parameter is done (see \cite{bib:nikiforov-uvarov-suslov} for classical Charlier polynomials orthogonal with respect to $\rho(k)$ and \cite{bib:arvesu-coussement-vanassche} for multiple Charlier polynomials).

From equations \eqref{eq:particular_measures}, \eqref{eq:particular_moments}, \eqref{eq:mass_points}, and \eqref{eq:ch_weight} we obtain a particular system of discrete measures $\{ \mu_0, \ldots, \mu_{r-1}\}$ supported on the $r$-star \eqref{r-star} with the following moments
\begin{equation}
m_{j}^{(\ell)} = \sum_{k=0}^\infty b_{k,j}^{(\ell)},\quad
b_{k,j}^{(\ell)}=k^j \dfrac{a_{\ell}^k}{\Gamma(k + 1)},\quad \ell=0,\ldots,r-1,\quad j\in \mathbb{N}.
\label{eq:ch_moments}
\end{equation}
For any fixed $j \in \mathbb{N}$ and $\ell = 0,...,r-1$, one has
\begin{equation*}
\left\lvert
b_{k+1,j}^{(\ell)}/b_{k,j}^{(\ell)}\right\rvert = \left\lvert a_\ell \right\rvert
\left\lvert  (k+1)^{j-1}/k^j \right\rvert,
\end{equation*}
which converges to zero as $k \rightarrow \infty$. From the Ratio Test follows the absolute convergence of the series, then \eqref{eq:ch_moments} converges and all moments exist. Despite this fact, due to the periodicity of the complex exponential function we need that
$$
\ln (a_j/a_\ell)\not=2\pi i k,\quad k\in\mathbb{Z},\quad \mbox{for}\quad j\not=\ell,\quad j,\ell=0,\ldots,r-1,
$$
in order the matrix $\mathcal{M}$ derived from the linear system \eqref{eq:type2_original} has full rank (see expressions \eqref{matrix} and \eqref{Pn_matrix_eq}). Thus, $\{\mu_\ell \}_{\ell=0}^{r-1}$ is a perfect system of measures.

\subsection{Meixner--Angelesco system of measures}
Consider the mass points \eqref{eq:mass_points}, the functions $\rho_{j}(z)$ in \eqref{eq:ch_weight}, in which the vector $\vec{a}$ is replaced by $\vec{c}=(c_0,\ldots,c_{r-1})\in \mathbb{C}^r$ with non-zero complex entries, and the weight functions
\begin{equation}\label{eq:m1_weight_function}
\rho_{j}^\beta (z) =\Gamma(z^r+\beta)\rho_j (z),\quad -\pi<\arg c_j\leq \pi,\quad j=0,\ldots, r-1,
\end{equation}
where the complex parameter $\beta$ is such that $(z^{r}_{j,k}+\beta)\in\mathbb{C}\setminus\mathbb{Z}_{-}$. In the expression \eqref{eq:m1_weight_function} a complex replacement of the real parameters involved in the negative binomial distribution (Pascal distribution) $\rho^{\beta}(k)=(\beta)_{k}c^{k}/k!$, on the non-negative integers $k=0,1,2\ldots$ is carried out (see \cite{bib:nikiforov-uvarov-suslov} and \cite{bib:arvesu-coussement-vanassche} for the real weight functions in Meixner polynomials and multiple Meixner polynomials, respectively).

From \eqref{eq:particular_measures}, \eqref{eq:mass_points}, and \eqref{eq:m1_weight_function} we obtain a particular system of discrete measures of Meixner of the first kind $\{ \mu_{0}^{\beta}, \ldots, \mu_{r-1}^{\beta}\}$ supported on the $r$-star \eqref{r-star} with the following moments
\begin{equation}
m_{j}^{(\ell)} = \sum_{k=0}^\infty b_{k,j}^{(\ell)},\quad
b_{k,j}^{(\ell)}
=k^j \dfrac{\Gamma(k + \beta)}{\Gamma(k + 1)}c_{\ell}^k,,\quad \ell=0,\ldots,r-1,\quad j\in \mathbb{N}.
\label{eq:m1_moments}
\end{equation}
For any fixed $j \in \mathbb{N}$ and $\ell = 0,...,r-1$, one gets
\begin{equation*}
\left\lvert
b_{k+1,j}^{(\ell)}/b_{k,j}^{(\ell)}\right\rvert = \left\lvert c_{\ell} \right\rvert \left\lvert \dfrac{(k+\beta)(k+1)^{j-1}}{k^j} \right\rvert \longrightarrow \left\lvert c_{\ell} \right\rvert, \ \ \text{as} \ \ k \rightarrow \infty.
\end{equation*}

Taking $\left\lvert c_{\ell} \right\rvert\leq c < 1$ for all $\ell = 0, \ldots, r-1$, from the Ratio Test follows the absolute convergence of the series. Then, \eqref{eq:m1_moments} converges and all moments exist. Due to the periodicity of the complex exponential function we need that
$$
\ln (c_j/c_\ell)\not=2\pi i k,\quad k\in\mathbb{Z},\quad \mbox{for}\quad j\not=\ell,\quad j,\ell=0,\ldots,r-1,
$$
in order the matrix $\mathcal{M}$ derived from the linear system \eqref{eq:type2_original} has full rank (see expressions \eqref{matrix} and \eqref{Pn_matrix_eq}). Thus, the system of discrete measures of Meixner of the first kind is a perfect system of measures.

\section{Type II Charlier--Angelesco polynomials}\label{sec:ch2_examples}

The type II monic Charlier--Angelesco polynomial
$C^{\vec{a}}_{\vec{n}}$ for the multi-index
$\vec{n}\in\mathbb{Z}_{+}^{r}$ on the $r$-star and weight functions
\eqref{eq:ch_weight} is the polynomial of degree
$\lvert\vec{n}\rvert=n_0+\cdots+n_{r-1}$ in $z^r$ defined by the
orthogonality relations
\begin{equation} \label{eq:ch_orth_cond_1}
\displaystyle{\int\limits_{\Omega_\ell}} C^{\vec{a}}_{\vec{n}}(z)
\left( -z^r \right)_j d\mu_{\ell}=0, \quad j=0,\ldots,
n_{\ell}-1,\quad \ell=0,\ldots,r-1,
\end{equation}
or, equivalently,
\begin{align}\label{eq:ch_orth_cond_2}
\displaystyle\sum_{k=0}^{\infty} C^{\vec{a}}_{\vec{n}}(z_{\ell,k})
\Pochh{-z_{\ell,k}^r}{j} \rho_{\ell}(z_{\ell,k}) &= 0,\quad 0\leq j
\leq n_{\ell}-1,\quad 0\leq \ell  \leq r-1.
\end{align}

The orthogonality relations \eqref{eq:ch_orth_cond_1} and
\eqref{eq:ch_orth_cond_2} are a particular situation of
\eqref{eq:type2_original} and \eqref{eq:discrete_orth_conditions},
with the $\omega$-symmetric weight functions \eqref{eq:ch_weight}.

\subsection{Raising relation and Rodrigues-type formula}\
Now we address the \textit{Rodrigues-type formula} for the type II
Charlier--Angelesco polynomials; however, some preliminary results
are needed.

\begin{lemma} \label{lemma:ch_raising_operator}
The type II Charlier--Angelesco polynomials verify the following
raising relation
\begin{equation} \label{eq:ch_raising}
C^{\vec{a}}_{\vec{n}+\vec{e}_{\ell}} (z) = \dfrac{-a_{\ell}}{\rho_{\ell}(z)}\nabla\left[ C^{\vec{a}}_{\vec{n}}(z) \rho_{\ell}(z) \right]=\Phi_{\ell}C^{\vec{a}}_{\vec{n}}(z),
\end{equation}
where the difference operator $\Phi_{\ell} = -a_{\ell} \left(\dfrac{1}{\rho_{\ell}(z)}\nabla \rho_{\ell}(z)\right)$ is acting on the polynomial $C^{\vec{a}}_{\vec{n}}$ and the function $\rho_{\ell}(z)$ is given in \eqref{eq:ch_weight}.
\end{lemma}
The operator $\Phi_{\ell}$ is called \textit{raising operator}
because the ${\ell}$th component of the multi-index $\vec{n}$ is
increased by $1$. In general, for any power $k=0,1,\ldots$,
\begin{align*}
\Phi_{\ell}z^k&=-a_{\ell} \left(\dfrac{\Gamma(z^r+1)}{a_{\ell}^{z^r}}
\nabla \dfrac{a_{\ell}^{z^r}}{\Gamma(z^r+1)}\right)z^k\\
&=-a_{\ell}
\dfrac{\Gamma(z^r+1)}{a_{\ell}^{z^r}}\left(\dfrac{a_{\ell}^{z^r}z^k}{\Gamma(z^r+1)}-
\dfrac{a_{\ell}^{z^r-1}(z^k-1)}{\Gamma(z^r)}\right)=z^{r}(z^{k}-1)-a_{\ell}z^k,
\end{align*}
where $\deg\Phi_{\ell}z^k=r+k$, and for $k=r\lvert\vec{n}\rvert$
the resulting polynomial has degree $r(\lvert\vec{n}\rvert+1)$.

\begin{proof} Since $\nabla(z^r+1)_{j+1}$, which is equal to $(j+1)(z^r+1)_{j}$, is a polynomial of degree $j$ in $z^r$ one can rewrite the orthogonality relations \eqref{eq:ch_orth_cond_2} as follows
\begin{equation}\label{orth_transition_ch}
\displaystyle{\sum^\infty_{k=0}} C^{\vec{a}}_{\vec{n}}\left(z_{{\ell},k}\right) \nabla\left(-z_{{\ell},k}^r +1\right)_{j+1} \rho_{\ell} \left(z_{{\ell},k} \right)=0,\quad 0\leq j  \leq n_{\ell}-1,\quad 0\leq \ell  \leq r-1.
\end{equation}
Considering the explicit expression of the mass points \eqref{eq:mass_points} along the rays, one gets $z_{{\ell},k}^r = k$, for ${\ell} = 0, \ldots, r-1$, hence the orthogonality relations \eqref{orth_transition_ch} become
\begin{equation}\label{orth_new_ch}
\displaystyle\sum^\infty_{k=0} C^{\vec{a}}_{\vec{n}}(k)
\nabla\left(-k +1\right)_{j+1} \rho_{\ell,k}=0,\quad 0\leq j  \leq n_{\ell}-1,\quad 0\leq \ell  \leq r-1,
\end{equation}
where $\rho_{\ell,k}=a_{\ell}^{k}/\Gamma(k+1)$. Using summation by
parts and the relations
\begin{align}
\rho_{\ell,-1} = a_{\ell}^{-1}/\Gamma\left(0\right) =0,\quad \lim_{k
\rightarrow \infty} C^{\vec{a}}_{\vec{n}}\left(k \right) \left( - k
\right)_{j+1} \rho_{\ell,k}=0,\label{lim_moment_term}
\end{align}
where \eqref{lim_moment_term} follows from the nth-Term Test for the
absolutely convergent series  \eqref{eq:ch_moments}, one has that
for any two polynomials $P(k)$ and $Q(k)$
\begin{equation}\label{eq_operatos_D_nabla}
\sum_{k=0}^{\infty}Q(k)\rho_{\ell,k}\Delta P(k)=-\sum_{k=0}^{\infty}P(k)\nabla\big(Q(k)\rho_{\ell,k}\big),
\end{equation}
holds. Thus, from \eqref{relation_delta_nabla} and
\eqref{eq_operatos_D_nabla} by  taking $P(k)=\left(-k \right)_{j+1}$
and $Q(k)=C^{\vec{a}}_{\vec{n}}(k)$, the relations
\eqref{orth_new_ch} became
$$
0=\sum^\infty_{k=0} C^{\vec{a}}_{\vec{n}}(k) \nabla\left(-k
+1\right)_{j+1} \rho_{\ell,k} =\sum_{k=0}^{\infty}
C^{\vec{a}}_{\vec{n}}(k)\rho_{\ell,k}\Delta \left(-k\right)_{j+1}
=-\sum_{k=0}^{\infty}\left(-k\right)_{j+1}
\nabla\left(C^{\vec{a}}_{\vec{n}}(k) \rho_{\ell,k}\right).
$$
Note that
\begin{equation}\label{eq4recurrence}
\nabla\left(C^{\vec{a}}_{\vec{n}}(k)
\rho_{\ell,k}\right)=\mathcal{Q}_{\lvert\vec{n}\rvert+1}(k)\rho_{\ell,k},
\end{equation}
where $\mathcal{Q}_{\lvert\vec{n}\rvert+1}(k)= C^{\vec{a}}_{\vec{n}}  \left(k \right) - a_{\ell}^{-1} k  C^{\vec{a}}_{\vec{n}} \left(k -1\right)$, which is a polynomial of degree exactly $\lvert\vec{n}\rvert+1$. Therefore, from the uniqueness of the monic multiple orthogonal polynomials derived from the orthogonality relations \eqref{orth_transition_ch} one gets $\mathcal{Q}_{\lvert\vec{n}\rvert+1}(z)=- a_{\ell}^{-1}C^{\vec{a}}_{\vec{n}+\vec{e}_{\ell}} (z)$, that is,
$$
\nabla\left(C^{\vec{a}}_{\vec{n}}(z)
\rho_{\ell}(z)\right)=- a_{\ell}^{-1}C^{\vec{a}}_{\vec{n}+\vec{e}_{\ell}} (z)\rho_{\ell}(z),
$$
which is equivalent to \eqref{eq:ch_raising}.
\end{proof}

Observe that from \eqref{eq:ch_raising} and \eqref{eq4recurrence}
one has the  following recurrence relation involving the type II
Charlier--Angelesco polynomials in the variable $z^r$,
$
C^{\vec{a}}_{\vec{n}+\vec{e}_{\ell}} (z)=z^{r}  C^{\vec{a}}_{\vec{n}} \left(z -1\right)- a_{\ell}C^{\vec{a}}_{\vec{n}}  \left(z \right).
$

\begin{lemma}\label{important_commutative_property} Let $a_{\ell}$ and
$a_{j}$ $(j,\ell=0,\ldots,r-1)$ be two entries of the vector $\vec{a}$ with
components $(a_0,\ldots,a_{r-1})\in \mathbb{C}^r$ given in \eqref{eq:ch_weight}. Then,
the following commutative relation
\begin{equation}\label{commutation_important_1}
\left(\dfrac{1}{a_{\ell}^{z^r}}\nabla a_{\ell}^{z^r}
\right)\left(\dfrac{1}{a_{j}^{z^r}}\nabla a_{j}^{z^r}
\right)=\left(\dfrac{1}{a_{j}^{z^r}}\nabla a_{j}^{z^r} \right)
\left(\dfrac{1}{a_{\ell}^{z^r}}\nabla a_{\ell}^{z^r} \right)
\end{equation}
holds.
\end{lemma}
\begin{proof} For convenience, let's take a continuous function $f(z)$.
A straightforward calculation leads to the following relation
\begin{align}
\left(\dfrac{1}{ a_{\ell}^{z^{r}}}\nabla a_{\ell}^{z^{r}}
\right) \cdot \left(\dfrac{1}{ a_{j}^{z^{r}}}\nabla a_{j}^{z^{r}}\right) f(z)
&= a_{\ell}^{-z^{r}}\nabla \left(
       a_{\ell}^{z^{r}} \cdot
       \dfrac{1}{ a_{j}^{z^{r}} }\nabla \left[a_{j}^{z^{r}} f(z)\right]\right)\notag\\
&= a_{\ell}^{-z^{r}}\nabla \left(a_{\ell}^{z^{r}}f(z)- a_{j}^{-1}a_{\ell}^{z^{r}} f(z-1)\right)\notag\\
&=f(z)+\dfrac{a_{\ell}+a_{j}}{a_{\ell}a_{j}}f(z-1)+\dfrac{1}{a_{\ell}a_{j}}f(z-2).\label{eq:ch_symmetric}
\end{align}
Since the last expression \eqref{eq:ch_symmetric} is invariant
under the changes $\ell$ by $j$ and vice versa, one gets
$$
f(z)+\dfrac{a_{j}+a_{\ell}}{a_{j}a_{\ell}}f(z-1)+\dfrac{1}{a_{j}a_{\ell}}f(z-2)=\left(\dfrac{1}{ a_{j}^{z^{r}}}\nabla a_{j}^{z^{r}}\right)\cdot\left(\dfrac{1}{ a_{\ell}^{z^{r}}}\nabla a_{\ell}^{z^{r}}
\right)f(z),
$$
which proves the relation \eqref{commutation_important_1}.
\end{proof}

Notice that for the product of $k$ equal terms one has
$$
\left(\dfrac{1}{a_{\ell}^{z^r}}\nabla a_{\ell}^{z^r} \right)\cdots\left(\dfrac{1}{a_{\ell}^{z^r}}\nabla a_{\ell}^{z^r}\right)=
\left(\dfrac{1}{a_{\ell}^{z^r}}\nabla^{k} a_{\ell}^{z^r} \right).
$$
Thus, from Lemma \ref{important_commutative_property} the following commutative relation
\begin{equation}\label{commutation_important}
\left(\dfrac{1}{a_{\ell}^{z^r}}\nabla^{m} a_{\ell}^{z^r} \right)\left(\dfrac{1}{a_{j}^{z^r}}\nabla^{k} a_{j}^{z^r} \right)
=\left(\dfrac{1}{a_{j}^{z^r}}\nabla^{k} a_{j}^{z^r} \right)
\left(\dfrac{1}{a_{\ell}^{z^r}}\nabla^{m} a_{\ell}^{z^r} \right),\,\, k,m\in\mathbb{N},
\end{equation}
holds.

\begin{corollary} \label{prop:ch_nabla_conmutativity}
Consider $\hat{\Phi}_\ell=-a^{-1}_{\ell}\Phi_\ell$ and the product
of raising operators $\hat{\Phi}_\ell \cdot \hat{\Phi}_j$ as a
composition $\left(\hat{\Phi}_\ell\circ \hat{\Phi}_j\right)$, i.e.,
\begin{equation*}
\left(\dfrac{1}{\rho_{\ell}(z)}\nabla \rho_{\ell}(z)\right) \left(\dfrac{1}{\rho_j(z)}\nabla \rho_j(z)\right)= \dfrac{1}{\rho_{\ell}(z)}\nabla \left[\rho_{\ell}(z) \left(\dfrac{1}{\rho_j(z)}\nabla \left(\rho_j(z)\right)\right)\right].
\end{equation*}
Then, the commutative property
\begin{equation}\label{eq_comm_phi}
\hat{\Phi}_\ell \cdot \hat{\Phi}_j = \hat{\Phi}_j \cdot
\hat{\Phi}_\ell
\end{equation}
holds.
\end{corollary}

\begin{proof} Using the explicit expressions for the weight functions  \eqref{eq:ch_weight} one gets
\begin{multline}\label{eq_4_commutation}
\hat{\Phi}_\ell \cdot \hat{\Phi}_j=\left( \dfrac{1}{\rho_{\ell}(z)}\nabla\rho_{\ell}(z)\right)\left( \dfrac{1}{\rho_j(z)}\nabla\rho_j(z)\right)
\\
= \dfrac{\Gamma(z^r+1)}{a_{\ell}^{z^r}} \nabla\left[ \dfrac{a_{\ell}^{z^r}}{\Gamma(z^r+1)}
\dfrac{\Gamma(z^r+1)}{a_j^{z^r}}\nabla \dfrac{a_j^{z^r}}{\Gamma(z^r+1)} \right] \\
= \Gamma(z^r+1) \left(  \dfrac{1}{a_{\ell}^{z^r}} \nabla  a_{\ell}^{z^r}\right)\left(  \dfrac{1}{a_j^{z^r}} \nabla  a_j^{z^r}\right)\left(\dfrac{1}{\Gamma(z^r+1)}\right).
\end{multline}
Since the operators $\left(  \dfrac{1}{a_{\ell}^{z^r}} \nabla  a_{\ell}^{z^r}\right)$ and $\left(  \dfrac{1}{a_j^{z^r}} \nabla  a_j^{z^r}\right)$ in \eqref{eq_4_commutation} are commuting (see Lemma \ref{important_commutative_property}) the relation
$\hat{\Phi}_\ell \cdot \hat{\Phi}_j = \hat{\Phi}_j \cdot \hat{\Phi}_\ell$ holds.
\end{proof}

In particular for the Charlier--Angelesco case, the equation
\eqref{eq_comm_phi} can be derived from the following relation
\begin{multline}\label{eq:ch_symmetry_weights_ralation}
\dfrac{\rho_{\ell}(z-1)}{\rho_{\ell}(z)}\dfrac{\rho_j(z-2)}{\rho_j(z-1)} = \dfrac{a_{\ell}^{z^r-1}}{\Gamma(z^r)} \dfrac{\Gamma(z^r+1)}{a_{\ell}^{z^r}} \dfrac{a_j^{z^r-2}}{\Gamma(z^r-1)}\dfrac{\Gamma(z^r)}{a_j^{z^r-1}} \\
    = \dfrac{a_{\ell}^{z^r-2}}{\Gamma(z^r-1)} \dfrac{\Gamma(z^r)}{a_{\ell}^{z^r-1}} \dfrac{a_j^{z^r-1}}{\Gamma(z^r)}\dfrac{\Gamma(z^r+1)}{a_j^{z^r}}= \dfrac{\rho_{\ell}(z-2)}{\rho_{\ell}(z-1)}\dfrac{\rho_j(z-1)}{\rho_j(z)}.
\end{multline}
Indeed,
\begin{align*}
\hat{\Phi}_\ell \cdot \hat{\Phi}_j (f) =
\left(\dfrac{1}{\rho_{\ell}(z)}\nabla \rho_{\ell}(z)\right) \cdot
\left(\dfrac{1}{\rho_j(z)}\nabla \rho_j(z)\right) f(z) =
\dfrac{1}{\rho_{\ell}(z)}\nabla \left[\rho_{\ell}(z) \cdot
\dfrac{1}{\rho_j(z)}\nabla \left(\rho_j(z) f(z)\right)\right].
\end{align*}
Hence
\begin{multline}\label{eq:ch_ustify_change}
\hat{\Phi}_\ell \cdot \hat{\Phi}_j (f)= \dfrac{1}{\rho_{\ell}(z)}\nabla
\left[ \dfrac{\rho_{\ell} (z)}{\rho_j(z)} \left(\rho_j(z)f(z)-\rho_j(z-1)f(z-1)\right) \right] \\
= f(z) - \left( \dfrac{\rho_{\ell}(z-1)}{\rho_{\ell}(z)} + \dfrac{\rho_j(z-1)}{\rho_j(z)}\right)f(z-1)
  +   \dfrac{\rho_{\ell}(z-2)}{\rho_{\ell}(z-1)}\dfrac{\rho_j(z-1)}{\rho_j(z)}f(z-2).
\end{multline}
Now using \eqref{eq:ch_symmetry_weights_ralation}, the equation \eqref{eq:ch_ustify_change} becomes
$\hat{\Phi}_\ell \cdot \hat{\Phi}_j (f)=\hat{\Phi}_j \cdot \hat{\Phi}_{\ell} (f)$.

\begin{theorem}\label{theorem:ch_rodrigues}
The type II Charlier--Angelesco polynomials are given by the
Rodrigues formula
\begin{align}\label{eq:ch_rodrigues_util}
C^{\vec{a}}_{\vec{n}} (z) &= \left(\prod_{\ell=0}^{r-1} \left(-a_{\ell}\right)^{n_{\ell}}
\right) \prod_{j=0}^{r-1} \left( \dfrac{1}{\rho_j(z^r)}\nabla^{n_j}\rho_j(z^r) \right)
\notag\\
&=  \left(\prod_{\ell=0}^{r-1}
\left(-a_{\ell}\right)^{n_{\ell}}\right) \Gamma(z^r+1)
\prod_{j=0}^{r-1} \left( \dfrac{1}{a_j^{z^r}}\nabla^{n_j}a_j^{z^r}
\right) \left(\dfrac{1}{\Gamma(z^r+1)}\right).
\end{align}
\end{theorem}

\begin{proof} In \eqref{eq:ch_raising} consider $\vec{n} = \vec{0}$ and
the monic polynomial $C^{\vec{a}}_{\vec{0}} (z)=1$. Hence,
\begin{align}\label{initial_exp_rod}
C^{\vec{a}}_{\vec{e}_\ell} (z) &= -a_{\ell} \dfrac{1}{\rho_{\ell}(z)}\nabla \rho_{\ell}(z)=
\Phi_{\ell}\cdot 1.
\end{align}
Multiplying \eqref{initial_exp_rod} by the product of
$(n_{\ell}-1)$-raising  operators $\Phi_{\ell}\cdots\Phi_{\ell}$,
one gets
\begin{align}\label{eq:ch_rodrigues_proof_one_ray}
C^{\vec{a}}_{n_{\ell}\vec{e}_\ell} (z) =  \Phi_{\ell}\cdots
\Phi_{\ell}\cdot 1= \left(-a_{\ell}\right)^{n_{\ell}} \left(
\dfrac{1}{\rho_{\ell}(z)}\nabla^{n_{\ell}} \rho_{\ell}(z)\right),
\end{align}
which follows from Lemma \ref{lemma:ch_raising_operator} and the
orthogonality conditions \eqref{eq:ch_orth_cond_2} for the
multi-index $n_{\ell}\vec{e}_\ell$. Similarly, for any other ray
$j\not=\ell$ of the $r$-star one has
\begin{multline*}
C^{\vec{a}}_{n_{j}\vec{e}_j+n_{\ell}\vec{e}_\ell} (z)=\underbrace{\Phi_{j}\cdots
\Phi_{j}}_{n_{j}-\mbox{times}}C^{\vec{a}}_{n_{\ell}\vec{e}_\ell} (z) 
= \left(-a_{j}\right)^{n_{j}}\left(-a_{\ell}\right)^{n_{\ell}}
\left( \dfrac{1}{\rho_{j}(z)}\nabla^{n_{j}} \rho_{j}(z)\right)
 \left( \dfrac{1}{\rho_{\ell}(z)}\nabla^{n_{\ell}} \rho_{\ell}(z)\right).
\end{multline*}
Thus, formula  \eqref{eq:ch_rodrigues_util} can be proved by induction. Indeed, assuming that
\begin{equation*}
C^{\vec{a}}_{n_{0}\vec{e}_0+\cdots+n_{r-2}\vec{e}_{r-2}} (z)=
\left(\prod_{\ell=0}^{r-2} \left(-a_{\ell}\right)^{n_{\ell}}\right) \prod_{j=0}^{r-2} \left(\dfrac{1}{\rho_j(z)}\nabla^{n_j}\rho_j(z) \right),
\end{equation*}
is valid and repeating the above process for the $(r-1)$-ray, one
gets the Rodrigues  formula  \eqref{eq:ch_rodrigues_util}. Note that
form Lemma \ref{prop:ch_nabla_conmutativity} the above inductive
method can be implemented in any order since the raising operators
are commuting.
\end{proof}

\subsection{Explicit expression and recurrence relation}

The explicit expression for the type II Charlier--Angelesco
polynomials will be derived from its Rodrigues formula
\eqref{eq:ch_rodrigues_util}. The next result is a corollary of
Theorem \ref{theorem:ch_rodrigues}.

\begin{corollary} The following explicit expression for type II Charlier--Angelesco
polynomials on the $r$-star \eqref{r-star} holds
\begin{equation} \label{eq:ch_explicit_expression}
C^{\vec{a}}_{\vec{n}} (z) = \sum_{k_0=0}^{n_{0}}\cdots
\sum_{k_{r-1}=0}^{n_{r-1}} \left(\prod_{\ell=0}^{r-1}
{n_{\ell}\choose{k_{\ell}}}\left(-a_{\ell}
\right)^{n_{\ell}-k_{\ell}}\right)(-1)^{\lvert\vec{k}\rvert}
\Pochh{-z^r}{\lvert\vec{k}\rvert}.
\end{equation}
\end{corollary}

\begin{proof} From the commutative property proved in Lemma \ref{prop:ch_nabla_conmutativity} and the Rodrigues formula \eqref{eq:ch_rodrigues_util} one has
\begin{align}\label{eq:ch_rodrigues_util_modif}
C^{\vec{a}}_{\vec{n}} (z) &=  \prod_{\ell=0}^{r-1} \left( \dfrac{\left(-a_{\ell}\right)^{n_{\ell}}}{\rho_{\ell}(z)}\nabla^{n_{\ell}}\rho_{\ell}(z) \right)=
\left(\prod_{\substack{\ell=0 \\ \ell\neq j}}^{r-1}\Phi_{\ell}^{n_{\ell}}\right) \Phi_{j}^{n_{j}}.
\end{align}
From \eqref{eq:ch_rodrigues_proof_one_ray} one has
\begin{align}\label{eq:ch_rodrigues_util_modif_1}
\Phi_{j}^{n_{j}}=C_{n_j\vec{e}_j}^{\vec{a}} (z)&= (-a_j)^{n_j} \Gamma(z^r+1) \left(\dfrac{1}
{a_j^{z^r}}\nabla^{n_j} a_j^{z^r}\right) \left(\dfrac{1}{\Gamma(z^r+1)}\right) \notag\\
&= \displaystyle\sum_{k_j=0}^{n_j} {n_j\choose{k_j}}(-a_j)^{n_j-k_j}\Pochh{z^r+1-{k_{j}} }{k_{j}} 
=   \displaystyle\sum_{k_j=0}^{n_j}
{n_j\choose{k_j}}(-a_j)^{n_j-k_j}(-1)^{k_{j}}\Pochh{-z^r}{k_{j}}.
\end{align}
In \eqref{eq:ch_rodrigues_util_modif_1} we use the relation $
\Pochh{z+1-n}{n} = (-1)^n \Pochh{-z}{n},$ where $ n \in \mathbb{N}.$
Therefore, equation \eqref{eq:ch_rodrigues_util_modif} becomes
\begin{align}\label{eq:ch_rodrigues_util_modif_2}
C^{\vec{a}}_{\vec{n}} (z) &=
\left(\prod_{\substack{\ell=0 \\ \ell\neq j}}^{r-1}\Phi_{\ell}^{n_{\ell}}\right)
\displaystyle\sum_{k_j=0}^{n_j} {n_j\choose{k_j}}(-a_j)^{n_j-k_j}(-1)^{k_{j}}\Pochh{-z^r}{k_{j}}.
\end{align}
Repeating the above process $(r-1)$-times in equation
\eqref{eq:ch_rodrigues_util_modif_2}  and taking into account that
the difference operators $\Phi_{\ell}^{n_{\ell}}$ can be taken in
any order since they are commuting operators, the formula
\eqref{eq:ch_explicit_expression} is obtained.
\end{proof}

The expression \eqref{eq:ch_explicit_expression} is given
explicitly as a linear combination of the basis elements
$\left\lbrace \Pochh{-z^r}{j} \right\rbrace_{j=0}^\infty$ formed
with Pochhammer symbols. In addition, one can easily check that the
resulting polynomials are monic.

Now, we focus our attention on the calculation of the explicit
expressions for the coefficients of the recurrence relation given in
Theorem \ref{theorem:recurrence_relation}. The next corollary is a
particular situation of \eqref{eq:recurrence_relation} for the set
of discrete measures $\mu_0, \ldots, \mu_{r-1}$ defined by the
weight functions \eqref{eq:ch_weight}.

\begin{corollary}\label{cor:ch_recurrence} The type II Charlier--Angelesco
polynomials verify the following recurrence relation
\begin{equation}\label{eq:ch_recurrence_relation}
\left(z^r -  b_{\vec{n},k}\right) C^{\vec{a}}_{\vec{n}}(z) - C_{\vec{n}+\vec{e}_k}^{\vec{a}}(z) = \sum_{j=0}^{r-1} d_{\vec{n},j}  C_{\vec{n}-\vec{e}_j}^{\vec{a}}(z),\quad k=0,\ldots,r-1,
\end{equation}
where $b_{\vec{n},k} = a_k + \lvert \vec{n}\rvert $ and $d_{\vec{n},j}=a_j n_j$.
\end{corollary}

\begin{proof} In Theorem \ref{theorem:recurrence_relation} was given $b_{\vec{n},k} =
\alpha_{\lvert \vec{n}\rvert -1}^{\vec{n}} - \alpha_{\lvert
\vec{n}\rvert }^{\vec{n}+\vec{e}_k}$.  This relation was derived
form the comparison of the coefficients of two near-neighbor
Charlier--Angelesco polynomials with multi-indices $\vec{n}$ and
$\vec{n}+\vec{e}_k$, written as a linear combination of the
monomials $(z^r)^m$, $0\leq m\leq \lvert\vec{n}\rvert+1$, i.e.,
\begin{equation}\label{eq:ch_formal_powers}
C^{\vec{a}}_{\vec{n}}(z) = \sum\limits_{m=0}^{\lvert\vec{n}\rvert}
\alpha_{m}^{\vec{n}} z^{rm}, \qquad
C_{\vec{n}+\vec{e}_k}^{\vec{a}}(z) =
\sum\limits_{m=0}^{\lvert\vec{n}\rvert+1}
\alpha_{m}^{\vec{n}+\vec{e}_k} z^{rm},
\end{equation}
where $\alpha_{\lvert \vec{n}\rvert }^{\vec{n}} = \alpha_{\lvert \vec{n}\rvert +1}^{\vec{n}+\vec{e}_k} = 1$.

Since the explicit expression for the type II Charlier--Angelesco
polynomials $C^{\vec{a}}_{\vec{n}}$ in
\eqref{eq:ch_explicit_expression} was given in terms of the
Pochhammer symbols one needs to implement a change of basis. Note
that
\begin{align}\label{eq:ch_pochhammer}
C^{\vec{a}}_{\vec{n}} (z) = (-1)^{\lvert \vec{n}\rvert }\Pochh{-z^r}{\lvert \vec{n}\rvert } + (-1)^{\lvert \vec{n}\rvert -1}\sum\limits_{j=0}^{r-1} (-a_j n_j) \Pochh{-z^r}{\lvert \vec{n}\rvert -1}+ \cdots
\end{align}
Now, using the relation (see \cite{bib:nikiforov-uvarov-suslov}, p. 53)
$
(-z)_{n}=(-1)^{n}\sum_{m=0}^{n}S_{n}^{(m)}z^{m},
$
where $S_{n}^{(m)}$ are the Stirling numbers of the first kind \cite{bib:abramowitz-stegun},
one can check that the first two coefficients in the expansion of $(-1)^n \Pochh{-z^r}{n}$ in terms of $\{(z^r)^n,(z^r)^{n-1},\ldots, 1\}$ are as follows
\begin{align}\label{eq:Pochh_second_coeff}
(-1)^n \Pochh{-z^r}{n} &=z^{rn} - {n\choose{2}}z^{rn-r} + \cdots
\end{align}
Similarly,
\begin{align*}
(-1)^{n+1} \Pochh{-z^r}{n+1}
&=z^{rn+r} - {n+1\choose{2}}z^{rn} + \cdots
\end{align*}
Therefore, using equation \eqref{eq:Pochh_second_coeff} in formula
\eqref{eq:ch_pochhammer} one gets the following explicit expressions
involving the main coefficients in \eqref{eq:ch_formal_powers}
\begin{align*}
C^{\vec{a}}_{\vec{n}} (z) &= z^{r\lvert \vec{n}\rvert } + \alpha_{\lvert \vec{n}\rvert -1}^{\vec{n}}z^{r\lvert \vec{n}\rvert -r} +\cdots +C^{\vec{a}}_{\vec{n}}(0),\\
C_{\vec{n}+\vec{e}_k}^{\vec{a}}(z) &= z^{r\lvert \vec{n}\rvert +r} + \alpha_{\lvert \vec{n}\rvert}^{\vec{n}+\vec{e}_k}z^{r\lvert \vec{n}\rvert } +\cdots +C_{\vec{n}+\vec{e}_k}^{\vec{a}}(0),
\end{align*}
where
\begin{align}
 \displaystyle\alpha_{\lvert \vec{n}\rvert -1}^{\vec{n}} &=-
 {\lvert \vec{n}\rvert\choose{2}} - \sum\limits_{j=0}^{r-1}a_j n_j ,\quad
 \alpha_{\lvert \vec{n}\rvert }^{\vec{n}+\vec{e}_k} =-
 {\lvert \vec{n}\rvert+1\choose{2}} - \sum\limits_{j=0}^{r-1}a_j n_j  -a_{k},\label{eq:ch_alpha_ns}\\
 C^{\vec{a}}_{\vec{n}}(0)&= \prod\limits_{\ell=0}^{r-1}\left(-a_\ell\right)^{n_\ell},\quad
 C^{\vec{a}}_{\vec{n}+\vec{e}_k}(0) = (-a_k)\prod\limits_{\ell=0}^{r-1}\left(-a_\ell\right)^{n_\ell}.\label{eq:ch_independent_term_ns}
\end{align}
Hence, from Theorem \ref{theorem:recurrence_relation}, equation
\eqref{eq:b_linear_combination1}, and \eqref{eq:ch_alpha_ns} $
b_{\vec{n},k} = \alpha_{\lvert \vec{n}\rvert -1}^{\vec{n}} -
\alpha_{\lvert \vec{n}\rvert }^{\vec{n}+\vec{e}_k} = a_k + \lvert
\vec{n}\rvert. $ For computing the coefficients $d_{\vec{n},j}$ we
evaluate the relation \eqref{eq:ch_recurrence_relation} at zero and
solve for $d_{\vec{n},j}$. Indeed,
\begin{equation}\label{eq:ch_recurrence}
-\left(a_k + \lvert \vec{n}\rvert  \right) C^{\vec{a}}_{\vec{n}}(0)  - C^{\vec{a}}_{\vec{n}+\vec{e}_k}(0) = \sum_{j=0}^{r-1} d_{\vec{n},j}  C^{\vec{a}}_{\vec{n}-\vec{e}_j}(0), \ \ \ k=0,\ldots,r-1,
\end{equation}
where
\begin{equation}\label{eq:ch_0_eval}
C^{\vec{a}}_{\vec{n}-\vec{e}_j}(0) = (-a_j)^{-1}\prod\limits_{\ell=0}^{r-1}\left(-a_\ell\right)^{n_\ell}.
\end{equation}
Substituting the independent terms given in
\eqref{eq:ch_independent_term_ns}  and \eqref{eq:ch_0_eval} in the
relation \eqref{eq:ch_recurrence} and dividing by
$\prod\limits_{\ell=0}^{r-1}\left(-a_\ell\right)^{n_\ell}$ one has
the following equation
\begin{align*}
\lvert \vec{n}\rvert &= \sum_{j=0}^{r-1} d_{\vec{n},j} (-a_j)^{-1}\quad
\Longleftrightarrow\quad \sum_{j=0}^{r-1} \left(n_j -d_{\vec{n},j} (-a_j)^{-1}\right)=0,
\end{align*}
which is valid for $d_{\vec{n},j} = n_j a_j$.
\end{proof}

\section{Type II Meixner--Angelesco polynomials of the first kind}\label{sec:meix2_examples}

The type II monic Meixner--Angelesco polynomial of the first kind
$M_{\vec{n}}^{\beta,\vec{c}}$ for the multi-index
$\vec{n}=(n_1,\ldots,n_r)$ on the $r$-star and weight functions
\eqref{eq:m1_weight_function} is the polynomial of degree
$\lvert\vec{n}\rvert$ in $z^r$ defined by the orthogonality
relations
\begin{equation} \label{eq:m1_orth_cond_1}
\displaystyle{\int\limits_{\Omega_\ell}}
M_{\vec{n}}^{\beta,\vec{c}}(z)  \left( -z^r \right)_j
d\mu_{\ell}^{\beta}=0, \quad j=0,\ldots, n_{\ell}-1,\quad
\ell=0,\ldots,r-1,
\end{equation}
or, equivalently,
\begin{align}\label{eq:m1_orth_cond_2}
\displaystyle\sum_{k=0}^{\infty} M_{\vec{n}}^{\beta,\vec{c}}(z_{\ell,k}) \Pochh{-z_{\ell,k}^r}{j} \rho_{\ell}^{\beta}(z_{\ell,k}) &= 0,\quad 0\leq j  \leq n_{\ell}-1,\quad 0\leq \ell  \leq r-1.
\end{align}

The orthogonality relations \eqref{eq:m1_orth_cond_1} and \eqref{eq:m1_orth_cond_2} are a particular situation of \eqref{eq:type2_original} and \eqref{eq:discrete_orth_conditions}, with the $\omega$-symmetric weight functions \eqref{eq:m1_weight_function}. Both systems \eqref{eq:m1_orth_cond_1} and \eqref{eq:m1_orth_cond_2} define $\lvert\vec{n}\rvert$ conditions for the $\lvert\vec{n}\rvert$-unknown coefficients of the monic polynomial $M_{\vec{n}}^{\beta,\vec{c}}$ of degree $\lvert\vec{n}\rvert$ in $z^r$.

\subsection{Raising relation and Rodrigues-type formula}

Now we address the \textit{raising operator} and
\textit{Rodrigues-type formula} for the type II Meixner--Angelesco
polynomials of the first kind.

Define
\begin{equation}\label{Psi_operator}
\Psi^{\beta}_\ell =\left( \dfrac{c_{\ell}}{c_{\ell}-1} \right) \dfrac{1}{\rho_{\ell}^{\beta-1}(z)}\nabla \rho_{\ell}^\beta(z),
\end{equation}
where $\rho_{\ell}^\beta(z)$ is given in \eqref{eq:m1_weight_function}. Indeed, for any power $k=0,1,\ldots$,
\begin{align}\label{Psi_power_m1}
\Psi^{\beta}_{\ell} z^k&=\dfrac{c_{\ell}}{c_{\ell}-1} \dfrac{\Gamma(z^r+1)}{c_{\ell}^{z^r}\Gamma(z^r+\beta-1)}
\left(\nabla \dfrac{c_{\ell}^{z^r}\Gamma(z^r+\beta)}{\Gamma(z^r+1)}z^k\right)\notag\\
&=\dfrac{c_{\ell}}{c_{\ell}-1} \dfrac{\Gamma(z^r+1)}{\Gamma(z^r+\beta-1)}\left(\dfrac{\Gamma(z^r+\beta)z^k}{\Gamma(z^r+1)}-
\dfrac{c_{\ell}^{-1}(z^k-1)\Gamma(z^r+\beta-1)}{\Gamma(z^r)}\right)\notag\\
&=z^{r+k}+\dfrac{c_{\ell}(\beta-1)}{c_{\ell}-1}z^{k}+\dfrac{1}{c_{\ell}-1}z^r,
\end{align}
where $\deg\Psi^{\beta}_{\ell}z^k=r+k$, and for $k=r\lvert\vec{n}\rvert$ the resulting polynomial has degree $r(\lvert\vec{n}\rvert+1)$.

\begin{lemma} \label{lemma:m1_raising_operator}
The type II Meixner--Angelesco polynomials of the first kind satisfy
the raising relation
\begin{equation} \label{eq:m1_raising}
M^{\beta-1,\vec{c}}_{\vec{n}+\vec{e}_{\ell}} \left(z\right) = \Psi^{\beta}_\ell M^{\beta,\vec{c}}_{\vec{n}} \left(z \right),
\end{equation}
where $\Psi^{\beta}_{\ell}$ is given in \eqref{Psi_operator}.
\end{lemma}

The operator $\Psi^{\beta}_\ell$ is called \textit{raising operator} because the ${\ell}$-th component of the multi-index $\vec{n}$ is increased by $1$.

\begin{proof} Replacing $(-z_{{\ell},k}^r)_{j}$ by $\nabla(-z_{{\ell},k}^r +1)_{j+1}$, $j=0,\ldots,n_{\ell}-1$, in the orthogonality relation \eqref{eq:m1_orth_cond_2} one has
\begin{equation}\label{orth_transition_m1}
\displaystyle{\sum^\infty_{k=0}} M^{\beta,\vec{c}}_{\vec{n}}\left(z_{{\ell},k}\right) \nabla\left(-z_{{\ell},k}^r +1\right)_{j+1} \rho^{\beta}_{\ell} \left(z_{{\ell},k} \right)=0,\quad 0\leq \ell  \leq r-1.
\end{equation}
Note that $\nabla(-z^r+1)_{j+1}$ is a polynomial of degree $j$ in $z^r$.

Considering the explicit expression of the mass points \eqref{eq:mass_points} along the rays, one gets $z_{{\ell},k}^r = k$, for ${\ell} = 0, \ldots, r-1$, hence the orthogonality relations \eqref{orth_transition_m1} become
\begin{equation}\label{orth_new_m1}
\displaystyle\sum^\infty_{k=0} M^{\beta,\vec{c}}_{\vec{n}}(k)
\nabla\left(-k +1\right)_{j+1} \rho^{\beta}_{\ell,k}=0,\quad 0\leq j  \leq n_{\ell}-1,\quad 0\leq \ell  \leq r-1,
\end{equation}
where $\rho^{\beta}_{\ell,k}=c_{\ell}^{k}\Gamma(\beta+k)/\Gamma(k+1)$.

In \eqref{orth_new_m1} we will use summation by parts and the relations
\begin{align}
\rho^{\beta}_{\ell,-1} =
c_{\ell}^{-1}\Gamma(\beta-1)/\Gamma\left(0\right) =0,\quad
 \lim_{k
\rightarrow \infty} M^{\beta,\vec{c}}_{\vec{n}}\left(k \right)
\left( - k  \right)_{j+1}
\rho^{\beta}_{\ell,k}=0,\label{lim_moment_term_m1}
\end{align}
where \eqref{lim_moment_term_m1} follows from the nth-Term Test for
the absolutely convergent series \eqref{eq:m1_moments}. Indeed, from
\eqref{relation_delta_nabla}  and \eqref{eq_operatos_D_nabla}, in
which we take $P(k)=\left(-k \right)_{j+1}$ and
$Q(k)=M^{\beta,\vec{c}}_{\vec{n}}(k)$, the relations
\eqref{orth_new_m1} became
\begin{align*}
0=\sum^\infty_{k=0} M^{\beta,\vec{c}}_{\vec{n}}(k) \nabla\left(-k
+1\right)_{j+1} \rho^{\beta}_{\ell,k} =\sum_{k=0}^{\infty}
M^{\beta,\vec{c}}_{\vec{n}}(k)\rho^{\beta}_{\ell,k}\Delta
\left(-k\right)_{j+1}=-\sum_{k=0}^{\infty}\left(-k\right)_{j+1}
\nabla\left(M^{\beta,\vec{c}}_{\vec{n}}(k)
\rho^{\beta}_{\ell,k}\right).
\end{align*}
Note that from \eqref{Psi_operator} and \eqref{Psi_power_m1}
\begin{equation}\label{eq4recurrence_m1}
\nabla\left(M^{\beta,\vec{c}}_{\vec{n}}(k)
\rho^{\beta}_{\ell,k}\right)=\mathcal{Q}_{\lvert\vec{n}\rvert+1}(k)\rho^{\beta-1}_{\ell,k},
\end{equation}
where $\mathcal{Q}_{\lvert\vec{n}\rvert+1}(k)=(k+\beta-1) M^{\beta,\vec{c}}_{\vec{n}}  \left(k \right) - \dfrac{k}{c_{\ell}}   M^{\beta,\vec{c}}_{\vec{n}} \left(k -1\right)$, which is a polynomial of degree exactly $\lvert\vec{n}\rvert+1$. Therefore, from the uniqueness of the monic multiple orthogonal polynomials derived from the orthogonality relations \eqref{orth_transition_m1} one gets $\mathcal{Q}_{\lvert\vec{n}\rvert+1}(z)=\left(\dfrac{c_{\ell}-1}{c_{\ell}}\right) M^{\beta-1,\vec{c}}_{\vec{n}+\vec{e}_{\ell}} \left( z \right)$, which is a polynomial of degree $\lvert\vec{n}\rvert+1$ in $z^r$. Accordingly, the expression \eqref{eq4recurrence_m1} becomes
$$
\nabla\left(M^{\beta,\vec{c}}_{\vec{n}}(z)
\rho^{\beta}_{\ell}(z)\right)=\left(\dfrac{c_{\ell}-1}{c_{\ell}}\right) M^{\beta-1,\vec{c}}_{\vec{n}+\vec{e}_{\ell}} (z)\rho^{\beta-1}_{\ell}(z).
$$
Finally, taking into account \eqref{Psi_operator} the raising relation \eqref{eq:m1_raising} holds.
\end{proof}

Observe that from \eqref{eq:m1_raising} and \eqref{eq4recurrence_m1}
one has the following recurrence relation involving the type II
Meixner--Angelesco polynomials in $z^r$,
$$
M^{\beta-1,\vec{c}}_{\vec{n}+\vec{e}_{\ell}} (z)=\left(
\dfrac{c_{\ell}} {c_{\ell}-1} \right)\left((z^r+\beta-1)
M^{\beta,\vec{c}}_{\vec{n}}  \left(z \right) - \dfrac{1}{c_{\ell}}
z^r  M^{\beta,\vec{c}}_{\vec{n}} \left(z -1\right)\right).
$$

In the sequel we will deal with compositions of raising operators involving the expression given in \eqref{Psi_operator}.

Define
\begin{align}\label{eq:m1_Psi_hat}
\hat{\Psi}^{\beta}_{\ell}&=\left(\dfrac{c_{\ell}-1}{c_{\ell}}\right)\Psi^{\beta}_{\ell}=\left(\dfrac{1}{\rho^{\beta-1}_{\ell}(z)}\nabla \rho^{\beta}_{\ell}(z)\right),\quad \ell=0,\ldots, r-1,
\end{align}
and consider the product of raising operators
$\hat{\Psi}^{\beta+1}_{\ell}  \cdots \hat{\Psi}^{\beta+k}_{\ell}$ as
a composition $\hat{\Psi}^{\beta+1}_{\ell}\circ\cdots\circ
\hat{\Psi}^{\beta+k}_{\ell}$, that is
\begin{align}\label{eq:m1_kth_times_operator_Psi}
\hat{\Psi}_{\ell}^{\beta+1}\cdots \hat{\Psi}_{\ell}^{\beta+k} =
\left( \dfrac{1}{\rho_{\ell}^{\beta}(z)}\nabla^{k}
\rho_{\ell}^{\beta+k} (z)\right)=\dfrac{\Gamma(z^r+1)
}{\Gamma(z^r+\beta) }\left(
\dfrac{1}{c_{\ell}^{z^r}}\nabla^{k}c_{\ell}^{z^r}
\right)\left(\dfrac{\Gamma(z^r+\beta+k
)}{\Gamma(z^r+1)}\right).
\end{align}

\begin{lemma} \label{prop:m1_nabla_conmutativity}
The following commutative property
\begin{multline}\label{operator_comm_m1}
\prod\limits_{s=0}^{k-1} \hat{\Psi}_{\ell}^{\beta+1 + s}\,
\prod\limits_{t=0}^{m-1} \hat{\Psi}_{j}^{\beta+1+k +
t}=\prod\limits_{t=0}^{m-1} \hat{\Psi}_{j}^{\beta+1 + t}\,
\prod\limits_{s=0}^{k-1} \hat{\Psi}_{\ell}^{\beta+1+m + s},\quad
0\leq j,\ell\leq r-1,\quad k,m\in\mathbb{N},
\end{multline}
holds.
\end{lemma}

\begin{proof} From \eqref{eq:m1_Psi_hat} the left-hand side of \eqref{operator_comm_m1} becomes
\begin{multline*}
\prod\limits_{s=0}^{k-1} \hat{\Psi}_{\ell}^{\beta+1 + s}\,
\prod\limits_{t=0}^{m-1} \hat{\Psi}_{j}^{\beta+1+k + t}= \left(
\dfrac{1}{\rho_{\ell}^{\beta}(z)}\nabla^{k} \rho_{\ell}^{\beta+k}
(z)\right)
\left( \dfrac{1}{\rho_{j}^{\beta+k}(z)}\nabla^{m} \rho_{j}^{\beta+k+m} (z)\right)\\
=\dfrac{\Gamma(z^r+1) }{\Gamma(z^r+\beta) }\left( \dfrac{1}{c_{\ell}^{z^r}}\nabla^{k}c_{\ell}^{z^r} \right)\left( \dfrac{1}{c_{j}^{z^r}}\nabla^{m}c_{j}^{z^r} \right)\left(\dfrac{\Gamma(z^r+\beta+k+m )}{\Gamma(z^r+1)}\right).
\end{multline*}
From \eqref{eq:m1_weight_function} the ratio
$\rho_{\ell}^{\beta+k}(z)/\rho_{j}^{\beta+k}(z)$ equals $1$.
Moreover, from Lemma \ref{important_commutative_property}, equation
\eqref{commutation_important}, the operators $\left(
\dfrac{1}{c_{\ell}^{z^r}}\nabla^{n_{\ell}}c_{\ell}^{z^r} \right)$
and $\left( \dfrac{1}{c_j^{z^r}}\nabla^{n_j}c_j^{z^r} \right)$ are
commuting. Therefore,
\begin{multline*}
\prod\limits_{s=0}^{k-1} \hat{\Psi}_{\ell}^{\beta+1 + s}\,
\prod\limits_{t=0}^{m-1} \hat{\Psi}_{j}^{\beta+1+k + t}
=\dfrac{\Gamma(z^r+1) }{\Gamma(z^r+\beta) }\left( \dfrac{1}{c_{j}^{z^r}}\nabla^{m}c_{j}^{z^r} \right)\left( \dfrac{1}{c_{\ell}^{z^r}}\nabla^{k}c_{\ell}^{z^r} \right)\left(\dfrac{\Gamma(z^r+\beta+k+m )}{\Gamma(z^r+1)}\right)\\
=\left( \dfrac{1}{\rho_{j}^{\beta}(z)}\nabla^{m} \rho_{j}^{\beta+m} (z)\right)
\left( \dfrac{1}{\rho_{\ell}^{\beta+m}(z)}\nabla^{k} \rho_{\ell}^{\beta+k+m} (z)\right)\\
=\prod\limits_{t=0}^{m-1} \hat{\Psi}_{j}^{\beta+1 + t}\,
\prod\limits_{s=0}^{k-1} \hat{\Psi}_{\ell}^{\beta+1+m + s},
\end{multline*}
which completes the proof.
\end{proof}

Let $\kappa=j_0,\ldots,j_{r-1}$ be any of the $r!$ permutations of $\{0,\ldots,r-1\}$ and $\vec{n}_{\kappa}=(n_{j_0},\ldots,n_{j_{r-1}})$, $\vec{n}=(n_0,\ldots,n_{r-1})$ the corresponding multi-indices, where $\lvert\vec{n}_{\kappa}\rvert=\lvert\vec{n}\rvert$. Thus, one can easily extend the above Lemma \ref{prop:m1_nabla_conmutativity} as follows
\begin{multline}\label{eq:m1_4_theorem_rodrigues}
\prod\limits_{t_{0}=0}^{n_{j_{0}}-1} \hat{\Psi}_{j_{0}}^{\beta+1 + t_{0}}
\prod\limits_{t_{1}=0}^{n_{j_{1}}-1} \hat{\Psi}_{j_{1}}^{\beta+1 +n_{j_{0}}+ t_{1}}\cdots
\prod\limits_{t_{r-1}=0}^{n_{j_{r-1}}-1} \hat{\Psi}_{j}^{\beta+1+\sum_{k=0}^{j_{r-2}}n_{k} + t_{r-1}}\\
=\dfrac{\Gamma(z^r+1) }{\Gamma(z^r+\beta) }
\left( \dfrac{1}{c_{j_{0}}^{z^r}}\nabla^{n_{j_{0}}}c_{j_{0}}^{z^r} \right)\cdots
\left( \dfrac{1}{c_{j_{r-1}}^{z^r}}\nabla^{n_{j_{r-1}}}c_{j_{r-1}}^{z^r} \right)
\dfrac{\Gamma(z^r+\beta + \lvert\vec{n}_{\kappa}\rvert)}{\Gamma(z^r+1)}\\
=\dfrac{\Gamma(z^r+1) }{\Gamma(z^r+\beta) }\prod_{j=0}^{r-1} \left( \dfrac{1}{c_{j}^{z^r}}
\nabla^{n_j}c_{j}^{z^r} \right)\dfrac{\Gamma(z^r+\beta + \lvert\vec{n}\rvert)}{\Gamma(z^r+1)}
=\prod\limits_{\ell=0}^{r-1}\left(\prod\limits_{j=0}^{n_{\ell}-1} \hat{\Psi}_{\ell}^{\beta+1 + j+\sum_{k=0}^{\ell}n_{k}-n_{\ell}}\right).
\end{multline}
Here we have arranged the product $\left(c_{j_{0}}^{-z^r}\nabla^{n_{j_{0}}}c_{j_{0}}^{z^r} \right)\cdots
\left( c_{j_{r-1}}^{-z^r}\nabla^{n_{j_{r-1}}}c_{j_{r-1}}^{z^r} \right)$ as $\displaystyle\prod_{j=0}^{r-1} \left( c_{j}^{-z^r}\nabla^{n_j}c_{j}^{z^r} \right)$ in accordance with the commutative property stated in Lemma \ref{important_commutative_property}.

\begin{theorem}\label{theorem:m1_rodrigues} The type II Meixner--Angelesco polynomials
of the first kind verify the Rodrigues-type formula
\begin{align}
M^{\beta, \vec{c}}_{\vec{n}} (z)&= \prod_{\ell=0}^{r-1} \left(\dfrac{c_{\ell}}{c_{\ell}-1}\right)^{n_\ell} \dfrac{\Gamma(z^r+1) }{\Gamma(z^r+\beta) }\prod_{j=0}^{r-1} \left( \dfrac{1}{c_j^{z^r}}\nabla^{n_j}c_j^{z^r} \right)\dfrac{\Gamma(z^r+\beta + \lvert\vec{n}\rvert)}{\Gamma(z^r+1)} \label{eq:m1_rodrigues_util}\\
&= \prod_{\ell=0}^{r-1} \left(\dfrac{c_{\ell}}{c_{\ell}-1}\right)^{n_\ell}\prod_{j=0}^{r-1} \left( \dfrac{1}{\rho_j^{\beta+\sum_{k=0}^j n_k -n_j}(z)}\nabla^{n_j} \rho_j^{\beta + \sum_{k=0}^j n_k}(z) \right).\notag
\end{align}
\end{theorem}

\begin{proof} In \eqref{eq:m1_raising} replace $\beta$ by $\beta+n_{\ell}$ (for any fixed ray $\ell\in\{0,\ldots, r-1\}$) and take the monic polynomial $M^{\beta+n_{\ell}, \vec{c}}_{\vec{0}} (z)=1$. Hence,
\begin{align}\label{m1_initial_exp_rod}
M^{\beta+n_{\ell}-1, \vec{c}}_{\vec{e}_{\ell}} (z) &= \Psi_{\ell}^{\beta+n_{\ell}} M^{\beta+n_{\ell}, \vec{c}}_{\vec{0}}= \left(\dfrac{c_{\ell}}{c_{\ell}-1} \right)\left(\dfrac{1}{\rho_{\ell}^{\beta+n_{\ell}-1}(z)}\nabla \rho_{\ell}^{\beta+n_{\ell}}(z)\right).
\end{align}
Multiplying equation \eqref{m1_initial_exp_rod} from the left by the
product of $(n_{\ell}-1)$-raising operators
$\Psi_{\ell}^{\beta+1}\cdots \Psi_{\ell}^{\beta+n_{\ell}-1}$ and
using \eqref{eq:m1_kth_times_operator_Psi}, one gets
\begin{align}\label{eq:m1_rodrigues_proof_one_ray}
M^{\beta, \vec{c}}_{n_{\ell}\vec{e}_{\ell}} (z)&=\prod\limits_{t=0}^{n_{\ell}-1} \Psi_{\ell}^{\beta+1 + t}=\left(\dfrac{c_{\ell}}{c_{\ell}-1}\right)^{n_{\ell}} \left( \dfrac{1}{\rho_{\ell}^{\beta}(z)}\nabla^{n_{\ell}} \rho_{\ell}^{\beta+n_{\ell}}(z)\right).
\end{align}
Here we have used Lemma \ref{lemma:m1_raising_operator} and the orthogonality conditions \eqref{eq:m1_orth_cond_2} for the multi-index $n_{\ell}\vec{e}_\ell$.

Similarly, take a ray $j\not=\ell$, $j=0,\ldots,r-1$,  and replace
in \eqref{eq:m1_rodrigues_proof_one_ray} $\beta$ by $\beta+n_j$,
then multiply \eqref{eq:m1_rodrigues_proof_one_ray} from the left by
the product of $n_{j}$-raising operators $\Psi_{j}^{\beta+1}\cdots
\Psi_{j}^{\beta+n_{j}}$ to obtain
\begin{align}\label{eq:m1_rodrigues_proof_second_ray}
M^{\beta, \vec{c}}_{n_{j}\vec{e}_{j}+n_{\ell}\vec{e}_{\ell}} (z)&=
\left(\prod\limits_{s=0}^{n_{j}-1} \Psi_{j}^{\beta+1 + s}\right)
M^{\beta+n_{j}, \vec{c}}_{n_{\ell}\vec{e}_{\ell}} (z)\notag\\
&=\left(\prod\limits_{s=0}^{n_{j}-1} \Psi_{j}^{\beta+1 + s}\right)
\left(\prod\limits_{t=0}^{n_{\ell}-1} \Psi_{\ell}^{\beta+1+n_{j} + t}\right)M^{\beta+n_{j}+n_{\ell}, \vec{c}}_{\vec{0}} (z)\notag\\
&=\dfrac{ c_{j}^{n_{j}} c_{\ell}^{n_{\ell}} }{ (c_{j}-1)^{n_{j}} (c_{\ell}-1)^{n_{\ell}}}
\left( \dfrac{1}{\rho_{j}^{\beta}(z)}\nabla^{n_{j}} \rho_{j}^{\beta+n_{j}}(z)\right)
 \left( \dfrac{1}{\rho_{\ell}^{\beta+n_{j}}(z)}\nabla^{n_{\ell}} \rho_{\ell}^{\beta+n_{j}+n_{\ell}}(z)\right),
\end{align}
where $M^{\beta+n_{j}+n_{\ell}, \vec{c}}_{\vec{0}} (z)=1$. Hence, using the explicit expression of the weight functions in \eqref{eq:m1_weight_function} one has
\begin{multline}\label{eq:m1_rodrigues_proof_second_ray_equiv}
M^{\beta, \vec{c}}_{n_{j}\vec{e}_{j}+n_{\ell}\vec{e}_{\ell}} (z)
=\dfrac{ c_{j}^{n_{j}} c_{\ell}^{n_{\ell}} }{ (c_{j}-1)^{n_{j}} (c_{\ell}-1)^{n_{\ell}}}\dfrac{\Gamma(z^r+1) }
{\Gamma(z^r+\beta) }\\
\times\left( \dfrac{1}{c_j^{z^r}}\nabla^{n_j}c_j^{z^r} \right)
\left( \dfrac{1}{c_{\ell}^{z^r}}\nabla^{n_{\ell}}c_{\ell}^{z^r}
\right)\dfrac{\Gamma(z^r+\beta + n_{j}+n_{\ell})}{\Gamma(z^r+1)}.
\end{multline}
In accordance with Lemmas \ref{important_commutative_property} and
\ref{prop:m1_nabla_conmutativity} the operators in
\eqref{eq:m1_rodrigues_proof_second_ray} and
\eqref{eq:m1_rodrigues_proof_second_ray_equiv} can be taken in any
order. Indeed, this is consistent with the fact that one could have
implemented the above procedure beginning with $j\in\{0,\ldots,
r-1\}$ and then proceeding with any $\ell\not=j$,
$\ell\in\{0,\ldots, r-1\}$, yielding the same result.

Repeating the above process for all other rays, taken in any order
and  arranged in accordance with \eqref{eq:m1_4_theorem_rodrigues}
one gets
\begin{align*}
M^{\beta, \vec{c}}_{\vec{n}} (z)&=\prod\limits_{\ell=0}^{r-1}\left(\prod\limits_{j=0}^{n_{\ell}-1}
\Psi_{\ell}^{\beta+1 + j+\sum_{k=0}^{\ell}n_{k}-n_{\ell}}\right)\\
&
= \prod_{\ell=0}^{r-1} \left(\dfrac{c_{\ell}}{c_{\ell}-1}\right)^{n_\ell}\prod_{j=0}^{r-1} \left( \dfrac{1}{\rho_j^{\beta+\sum_{k=0}^j n_k -n_j}(z)}\nabla^{n_j} \rho_j^{\beta + \sum_{k=0}^j n_k}(z) \right),
\end{align*}
which is equivalent to \eqref{eq:m1_rodrigues_util}.
\end{proof}

\subsection{Explicit expression and recurrence relation}\

Here we obtain the explicit expression for the type II
Meixner--Angelesco polynomials of the first kind from the Rodrigues
formula \eqref{eq:m1_rodrigues_util}. The next result is a
consequence of Theorem \ref{theorem:m1_rodrigues}.

\begin{corollary}\label{m1_quasi_final_corollary} The type II Meixner--Angelesco
polynomials of the first kind on the $r$-star are given by
\begin{equation} \label{eq:m1_explicit_expression}
M^{\beta, \vec{c}}_{\vec{n}}(z) =
\sum_{k_0=0}^{n_{0}}\cdots \sum_{k_{r-1}=0}^{n_{r-1}}\left(
\prod_{\ell=0}^{r-1} {n_{\ell}\choose{k_{\ell}}}
\dfrac{c_{\ell}^{n_{\ell}-k_{\ell}}}{\left( c_{\ell}-1 \right)^{n_{\ell}}} \right)
\Pochh{-z^r}{\lvert\vec{k}\rvert} \Pochh{z^r+\beta}{\lvert\vec{n}\rvert-\lvert\vec{k}\rvert},
\end{equation}
where $\lvert\vec{k}\rvert = k_0 + \ldots + k_{r-1}$.
\end{corollary}

\begin{proof} Consider the following calculation
\begin{multline}\label{eq:m1_initial_cal_explicit}
\left(\dfrac{1}{c_{r-1}^{z^r}}\nabla^{n_{r-1}} c_{r-1}^{z^r}\right)
\dfrac{\Gamma(z^r+\beta+\lvert\vec{n}\rvert)}{\Gamma(z^r+1)}\\
=\displaystyle{\sum_{k_{r-1}=0}^{n_{r-1}}
{n_{r-1}\choose{k_{r-1}}}}(-1)^{k_{r-1}}c_{r-1}^{-k_{r-1}}
\dfrac{\Gamma(z^r+\beta+\lvert\vec{n}\rvert-k_{r-1})}{\Gamma(z^r+1-k_{r-1})},
\end{multline}
in which we have used formula \eqref{nth_nabla}. Hence, base on the structure of the Rodrigues formula \eqref{eq:m1_rodrigues_util} and applying the above formula \eqref{eq:m1_initial_cal_explicit} adjusted to all other difference operators $\left(c_j^{-z^r}\nabla^{n_j}c_j^{z^r}\right)$, $j=0,\ldots,r-2$, involved in \eqref{eq:m1_rodrigues_util} one has
\begin{multline}\label{eq:m1_4_explicit}
\prod_{j=0}^{r-2} \left( \dfrac{1}{c_j^{z^r}}\nabla^{n_j}c_j^{z^r} \right)\left(\dfrac{1}{c_{r-1}^{z^r}}\nabla^{n_{r-1}} c_{r-1}^{z^r}\right) \dfrac{\Gamma(z^r+\beta+\lvert\vec{n}\rvert)}{\Gamma(z^r+1)}\\
=\sum_{k_0=0}^{n_{0}}\cdots \sum_{k_{r-1}=0}^{n_{r-1}}
\left(\prod_{\ell=0}^{r-1} {n_{\ell}\choose{k_{\ell}}}
c_{\ell}^{-k_{\ell}}\right) (-1)^{\lvert\vec{k}\rvert}\dfrac{\Gamma(z^r+\beta+\lvert\vec{n}\rvert-\lvert\vec{k}\rvert)}{\Gamma(z^r+1-\lvert\vec{k}\rvert)}.
\end{multline}
Note that the product of difference operators can be taken in any order.

Multiplying equation \eqref{eq:m1_4_explicit} from the left by $\prod_{\ell=0}^{r-1}\left(\dfrac{c_{\ell}}{c_{\ell}-1}\right)^{n_{\ell}}\dfrac{\Gamma(z^r+1)}{\Gamma(z^r+\beta)}$, one obtains
\begin{multline} \label{eq:m1_explicit_expression_1}
M^{\beta, \vec{c}}_{\vec{n}}(z) = \prod_{\ell=0}^{r-1}\left(\dfrac{c_{\ell}}
{c_{\ell}-1}\right)^{n_{\ell}}\dfrac{\Gamma(z^r+1)}{\Gamma(z^r+\beta)}\\
\times \sum_{k_0=0}^{n_{0}}\cdots \sum_{k_{r-1}=0}^{n_{r-1}}
\left(\prod_{\ell=0}^{r-1} {n_{\ell}\choose{k_{\ell}}}
c_{\ell}^{-k_{\ell}}\right) (-1)^{\lvert\vec{k}\rvert}
\dfrac{\Gamma(z^r+\beta+\lvert\vec{n}\rvert-\lvert\vec{k}\rvert)}{\Gamma(z^r+1-\lvert\vec{k}\rvert)}.
\end{multline}
In order to transform the right-hand side of equation \eqref{eq:m1_explicit_expression_1} we take into account the following relation
\begin{align*}
\dfrac{\Gamma(z^r+1) }{\Gamma(z^r+1-\lvert\vec{k}\rvert)}
\dfrac{\Gamma(z^r+\beta+\lvert\vec{n}\rvert-\lvert\vec{k}\rvert)}
{\Gamma(z^r+\beta) }
&=(-1)^{\lvert\vec{k}\rvert}\Pochh{-z^r}{\lvert\vec{k}\rvert} \Pochh{z^r+\beta}{\lvert\vec{n}\rvert-\lvert\vec{k}\rvert},
\end{align*}
for all $\lvert\vec{k}\rvert=0,\ldots,\lvert\vec{n}\rvert$, then
formula \eqref{eq:m1_explicit_expression_1} becomes
\eqref{eq:m1_explicit_expression}, i.e.,
\begin{align*}
M^{\beta, \vec{c}}_{\vec{n}}(z) &=
\sum_{k_0=0}^{n_{0}}\cdots \sum_{k_{r-1}=0}^{n_{r-1}}\left(
\prod_{\ell=0}^{r-1} {n_{\ell}\choose{k_{\ell}}}
\dfrac{c_{\ell}^{n_{\ell}-k_{\ell}}}{\left( c_{\ell}-1 \right)^{n_{\ell}}} \right)
\Pochh{-z^r}{\lvert\vec{k}\rvert} \Pochh{z^r+\beta}{\lvert\vec{n}\rvert-\lvert\vec{k}\rvert}.
\end{align*}
\end{proof}

In order to compute the coefficients of the recurrence relation for
the type II Meixner--Angelesco polynomials of the first kind given
in Theorem \ref{theorem:recurrence_relation} we use the explicit
expression obtained in Corollary \ref{m1_quasi_final_corollary}. The
procedure is similar to the one used in Corollary
\ref{cor:ch_recurrence}.

\begin{corollary} The type II Meixner--Angelesco polynomials of the first kind verify the
following recurrence relation
\begin{equation}\label{eq:m1_recurrence_relation}
\left(z^r - b_{\vec{n},\ell} \right) M^{\beta, \vec{c}}_{\vec{n}}(z) - M_{\vec{n}+\vec{e}_\ell}^{\beta, \vec{c}}(z) = \sum_{j=0}^{r-1} d_{\vec{n}, j} M_{\vec{n}-\vec{e}_j}^{\beta, \vec{c}}(z), \ \ \ \ell=0,\ldots,r-1,
\end{equation}
where $$b_{\vec{n},\ell} = \left(\lvert\vec{n}\rvert+\beta\right)\dfrac{c_\ell}{c_\ell -1 } + \sum\limits_{j=0}^{r-1} \dfrac{n_j}{c_j-1}\quad \text{and}\quad d_{\vec{n}, j} = c_j n_j \dfrac{\beta + \lvert\vec{n}\rvert - 1}{\left(c_j - 1\right)^2}.$$
\end{corollary}

\begin{proof} Base on Theorem \ref{theorem:recurrence_relation} we
calculate $b_{\vec{n},\ell} = \alpha_{\lvert \vec{n}\rvert
-1}^{\vec{n}} -  \alpha_{\lvert \vec{n}\rvert
}^{\vec{n}+\vec{e}_\ell}$ for the two near-neighbor
Meixner--Angelesco polynomials of the first kind
\begin{equation}\label{eq:m1_formal_powers}
M^{\beta, \vec{c}}_{\vec{n}}(z) = \sum\limits_{m=0}^{\lvert\vec{n}\rvert} \alpha_{m}^{\vec{n}} z^{rm}, \quad
M^{\beta, \vec{c}}_{\vec{n}+\vec{e}_\ell}(z)= \sum\limits_{m=0}^{\lvert\vec{n}\rvert+1} \alpha_{m}^{\vec{n}+\vec{e}_\ell} z^{rm}.
\end{equation}
Observe that the leading coefficient
\begin{align*}
\alpha_{\lvert\vec{n}\rvert}^{\lvert\vec{n}\rvert} =  \sum_{k_0=0}^{n_{0}}\cdots
\sum_{k_{r-1}=0}^{n_{r-1}}\prod_{j=0}^{r-1} {n_j\choose{k_j}}\dfrac{c_j^{n_j-k_j}}{\left( c_j-1
\right)^{n_i}}(-1)^{\lvert\vec{k}\rvert} 
    = \prod_{j=0}^{r-1} \dfrac{1}{\left( c_j-1 \right)^{n_j}} \sum_{k_j=0}^{n_{j}}{n_j\choose{k_j}}c_j^{n_j-k_j}(-1)^{k_j}= 1,
\end{align*}
and $\alpha_{\lvert \vec{n}\rvert }^{\vec{n}} = \alpha_{\lvert \vec{n}\rvert +1}^{\vec{n}+\vec{e}_\ell} = 1$.

From \eqref{eq:Pochh_second_coeff} and Newton's binomial formula we have
\begin{multline*}
\Pochh{z^r+\beta}{\lvert\vec{n}\rvert-\lvert\vec{k}\rvert} = (z^r+\beta)^{\lvert\vec{n}\rvert-\lvert\vec{k}\rvert } -{\lvert \vec{n}\rvert-\lvert\vec{k}\rvert\choose{2}}
(z^r+\beta)^{\lvert\vec{n}\rvert-\lvert\vec{k}\rvert -1}+\cdots\\
=z^{r(\lvert\vec{n}\rvert-\lvert\vec{k}\rvert)}+
\left({\lvert \vec{n}\rvert-\lvert\vec{k}\rvert\choose{1}}\beta+
{\lvert \vec{n}\rvert-\lvert\vec{k}\rvert\choose{2}}\right)
z^{r(\lvert\vec{n}\rvert-\lvert\vec{k}\rvert-1)}+\cdots
\end{multline*}
Hence,
\begin{multline}\label{eq:m1_Pochh_second_coeff}
\Pochh{-z^r}{\lvert\vec{k}\rvert} \Pochh{z^r+\beta}{\lvert\vec{n}\rvert-\lvert\vec{k}\rvert} =(-1)^{\lvert\vec{k}\rvert}z^{r\lvert\vec{n}\rvert}\\
+(-1)^{\lvert\vec{k}\rvert}z^{r\lvert\vec{n}\rvert-1}
\left({\lvert \vec{n}\rvert-\lvert\vec{k}\rvert\choose{1}}\beta+
{\lvert \vec{n}\rvert-\lvert\vec{k}\rvert\choose{2}}-
{\lvert \vec{k}\rvert\choose{2}}\right)+\cdots
\end{multline}
Therefore, using \eqref{eq:m1_explicit_expression} and
\eqref{eq:m1_Pochh_second_coeff} one obtains the following explicit
expressions  involving the main coefficients in
\eqref{eq:m1_formal_powers}
\begin{align*}
M^{\beta, \vec{c}}_{\vec{n}}(z) &= z^{r\lvert \vec{n}\rvert } + \alpha_{\lvert \vec{n}\rvert -1}^{\vec{n}}z^{r\lvert \vec{n}\rvert -r} +\cdots +M^{\beta, \vec{c}}_{\vec{n}}(0),\\
M^{\beta, \vec{c}}_{\vec{n}+\vec{e}_\ell}(z)&= z^{r\lvert \vec{n}\rvert +r} + \alpha_{\lvert \vec{n}\rvert}^{\vec{n}+\vec{e}_\ell}z^{r\lvert \vec{n}\rvert } +\cdots +M^{\beta, \vec{c}}_{\vec{n}+\vec{e}_\ell}(0),
\end{align*}
where
\begin{multline}\label{eq:m1_alpha_n-1}
\alpha_{\lvert \vec{n}\rvert -1}^{\vec{n}}=\left(\prod_{j=0}^{r-1} \dfrac{c_{j}^{n_{j}}}
{\left( c_{j}-1 \right)^{n_{j}}}\right)
\sum_{k_0=0}^{n_{0}}\cdots \sum_{k_{r-1}=0}^{n_{r-1}}
\left(\prod_{\ell=0}^{r-1} c^{-k_{\ell}}{n_{\ell}\choose{k_{\ell}}}\right)\\
\times
(-1)^{\lvert\vec{k}\rvert}
\left({\lvert \vec{n}\rvert-\lvert\vec{k}\rvert\choose{1}}\beta+
{\lvert \vec{n}\rvert-\lvert\vec{k}\rvert\choose{2}}-
{\lvert \vec{k}\rvert\choose{2}}\right),
\end{multline}
\begin{multline}\label{eq:m1_alpha_n}
\alpha_{\lvert \vec{n}\rvert }^{\vec{n}+\vec{e}_\ell}=\dfrac{c_{\ell}}
{\left( c_{\ell}-1 \right)}\left(\prod_{j=0}^{r-1} \dfrac{c_{j}^{n_{j}}}
{\left( c_{j}-1 \right)^{n_{j}}}\right)\\
\sum_{k_0=0}^{n_{0}} \cdots \sum_{k_{\ell}=0}^{n_{\ell}+1} \cdots \sum_{k_{r-1}=0}^{n_{r-1}}
\prod_{i=0}^{r-1} c^{-k_{i}}
(-1)^{\lvert\vec{k}\rvert}{n_{0}\choose{k_{0}}}\cdots {n_{\ell}+1\choose{k_{\ell}}}\cdots{n_{r-1}\choose{k_{r-1}}}\\
\times \left({\lvert \vec{n}\rvert+1-\lvert\vec{k}\rvert\choose{1}}\beta+
{\lvert \vec{n}\rvert+1-\lvert\vec{k}\rvert\choose{2}}-
{\lvert \vec{k}\rvert\choose{2}}\right),
\end{multline}
and
\begin{align}
M^{\beta, \vec{c}}_{\vec{n}}(0)&=(\beta)_{\lvert \vec{n}\rvert}\left(\prod_{j=0}^{r-1} \dfrac{c_{j}^{n_{j}}}
{\left( c_{j}-1 \right)^{n_{j}}}\right),\quad 
 M^{\beta, \vec{c}}_{\vec{n}+\vec{e}_\ell}(0)&=(\beta)_{\lvert \vec{n}\rvert+1}\dfrac{c_{\ell}}
{\left( c_{\ell}-1 \right)}\left(\prod_{j=0}^{r-1} \dfrac{c_{j}^{n_{j}}}
{\left( c_{j}-1 \right)^{n_{j}}}\right). \label{eq:m1_independent_term_ns}
\end{align}
Thus, from \eqref{eq:m1_alpha_n-1} and \eqref{eq:m1_alpha_n} one has
$$b_{\vec{n},\ell}= \alpha_{\lvert \vec{n}\rvert -1}^{\vec{n}} - \alpha_{\lvert \vec{n}\rvert }^{\vec{n}+\vec{e}_\ell} = -\left(\left(\lvert\vec{n}\rvert+\beta\right)\dfrac{c_\ell}{c_\ell -1 } + \sum\limits_{j=0}^{r-1} \dfrac{n_j}{c_j-1}\right).
$$

For computing the coefficients $d_{\vec{n},j}$ we evaluate the relation \eqref{eq:m1_recurrence_relation} at zero and solve for $d_{\vec{n},j}$. Indeed,
\begin{multline}\label{eq:m1_recurrence_relation_at_0}
  \left(\left(\lvert\vec{n}\rvert+\beta\right)\dfrac{c_\ell}{c_\ell -1 }
+ \sum\limits_{j=0}^{r-1} \dfrac{n_j}{c_j-1}\right)M^{\beta, \vec{c}}_{\vec{n}}(0) -
M_{\vec{n}+\vec{e}_\ell}^{\beta, \vec{c}}(0)\\
  = \sum_{j=0}^{r-1} d_{\vec{n}, j} M_{\vec{n}-\vec{e}_j}^{\beta, \vec{c}}(0), \ \  \ell=0,\ldots,r-1,
\end{multline}
where
\begin{equation}\label{eq:m1_0_eval}
M_{\vec{n}-\vec{e}_j}^{\beta, \vec{c}}(0) =(\beta)_{\lvert \vec{n}\rvert-1}\left(\dfrac{c_{j}-1}
{c_{j}}\right)\prod_{\ell=0}^{r-1} \dfrac{c_{\ell}^{n_{\ell}}}
{\left( c_{\ell}-1 \right)^{n_{\ell}}} .
\end{equation}
Substituting the independent terms given in \eqref{eq:m1_independent_term_ns} and \eqref{eq:m1_0_eval} in equation \eqref{eq:m1_recurrence_relation_at_0} and dividing both sides of the equation by $(\beta)_{\lvert \vec{n}\rvert-1}\prod_{\ell=0}^{r-1} c_{\ell}^{n_{\ell}}/
\left( c_{\ell}-1 \right)^{n_{\ell}}$ one gets the following equation
\begin{align*}
(\beta+\lvert \vec{n}\rvert-1)\sum_{j=0}^{r-1}\dfrac{n_j}{c_j-1} &=
\sum_{j=0}^{r-1} d_{\vec{n},j} \dfrac{c_{j}-1}{c_{j}}.
\end{align*}
Therefore,
$
d_{\vec{n},j} = \dfrac{n_j c_j}{(c_j-1)^2}(\beta+\lvert \vec{n}\rvert-1).
$
\end{proof}

\section{Conclusions}\label{conclusions}

Like the limit relations that involve the Charlier and Meixner
polynomials  (see formula (9.14.1) in
\cite{bib:koekoek-lesky-swarttouw}) we can obtain the type II
Charlier--Angelesco polynomials \eqref{eq:ch_explicit_expression} on
the $r$-star as limiting cases of type II Meixner--Angelesco
polynomials \eqref{eq:m1_explicit_expression}. For such a purpose we
will use the relation
\begin{equation}\label{eq:gamma_asymp}
\dfrac{\Gamma(w+a)}{\Gamma(w)}=w^{a}\left(1+O\left(w^{-1}\right)\right),\quad \lvert\arg w\rvert \leq \pi-\delta,\quad\delta>0.
\end{equation}
Furthermore, we replace $\vec{c}$ by $\vec{\gamma}=(a_{0}/(a_0+\beta),\ldots,a_{r-1}/(a_{r-1}+\beta))$; that is, in the equation \eqref{eq:m1_explicit_expression}
$c_\ell=(a_\ell+\beta)^{-1}a_{\ell}$, $\ell=0,\ldots,r-1$, and the coefficient
$$
\prod_{\ell=0}^{r-1} {n_{\ell}\choose{k_{\ell}}}
\dfrac{c_{\ell}^{n_{\ell}-k_{\ell}}}{\left( c_{\ell}-1 \right)^{n_{\ell}}} =(-1)^{\lvert\vec{n}\rvert}\beta^{\lvert\vec{k}\rvert-\lvert\vec{n}\rvert}\prod_{\ell=0}^{r-1} {n_{\ell}\choose{k_{\ell}}}a_{\ell}^{n_{\ell}-k_{\ell}}\left(1+\dfrac{a_{\ell}}{\beta}\right)^{k_{\ell}}.
$$
Indeed, equation \eqref{eq:m1_explicit_expression} becomes
\begin{multline} \label{eq:m1_explicit_expression_limit_charlier}
M^{\beta,\vec{\gamma}}_{\vec{n}}(z) =(-1)^{\lvert\vec{n}\rvert}\sum_{k_0=0}^{n_{0}}\cdots
\sum_{k_{r-1}=0}^{n_{r-1}}\left(\prod_{\ell=0}^{r-1} {n_{\ell}\choose{k_{\ell}}}a_{\ell}^{n_{\ell}-k_{\ell}}
\left(1+\dfrac{a_{\ell}}{\beta}\right)^{k_{\ell}}\right)\\
\times
\beta^{\lvert\vec{k}\rvert-\lvert\vec{n}\rvert}\Pochh{-z^r}{\lvert\vec{k}\rvert}
\Pochh{z^r+\beta}{\lvert\vec{n}\rvert-\lvert\vec{k}\rvert}.
\end{multline}
Now, taking into account \eqref{eq:gamma_asymp} one gets
\begin{align}\label{eq:m1_limit_1}
\lim_{\beta\to\infty}\beta^{\lvert\vec{k}\rvert-\lvert\vec{n}\rvert}\Pochh{z^r+\beta}{\lvert\vec{n}\rvert-\lvert\vec{k}\rvert}&=\lim_{\beta\to\infty}\beta^{\lvert\vec{k}\rvert-\lvert\vec{n}\rvert}
\dfrac{\Gamma(z^r+\beta+\lvert\vec{n}\rvert-\lvert\vec{k}\rvert)}
{\Gamma(z^r+\beta) }
=1.
\end{align}
Finally, from \eqref{eq:m1_explicit_expression_limit_charlier} and \eqref{eq:m1_limit_1} the following limit relation
\begin{multline*}
\lim_{\beta\to\infty}M^{\beta,\vec{\gamma}}_{\vec{n}}(z)=(-1)^{\lvert\vec{n}\rvert}\sum_{k_0=0}^{n_{0}}\cdots \sum_{k_{r-1}=0}^{n_{r-1}}\left(\prod_{\ell=0}^{r-1} {n_{\ell}\choose{k_{\ell}}}a_{\ell}^{n_{\ell}-k_{\ell}}\right) \Pochh{-z^r}{\lvert\vec{k}\rvert}\\
=\sum_{k_0=0}^{n_{0}}\cdots \sum_{k_{r-1}=0}^{n_{r-1}}\left(\prod_{\ell=0}^{r-1} {n_{\ell}\choose{k_{\ell}}}\left(-a_{\ell} \right)^{n_{\ell}-k_{\ell}}\right)(-1)^{\lvert\vec{k}\rvert} \Pochh{-z^r}{\lvert\vec{k}\rvert}=C^{\vec{a}}_{\vec{n}} (z),
\end{multline*}
holds.

In closing, we would like to point out a few comments, remarks, and open questions.

Some discrete multiple orthogonal polynomials for AT-systems of measures were studied in \cite{bib:arvesu-coussement-vanassche}. The main algebraic properties (raising operator, Rodrigues formula, recurrence relation for $r=2$) were addressed. In addition, it was indicated (see Conclusions p. 41) that it is not obvious to find Angelesco systems, where the discrete measures are supported on disjoint or touching intervals, which still have some raising operators, Rodrigues formula, and recurrence relations among other algebraic properties. The present paper has addressed that question. We have obtained the raising operators, the Rodrigues-type formulas, and the nearest neighbor recurrence relations, where the explicit expressions \eqref{eq:ch_explicit_expression} and \eqref{eq:m1_explicit_expression} played an important role. Note that these expressions represent polynomials in $z^r$, which differ from those of multiple Charlier polynomials and multiple Meixner polynomials of the first kind given in \cite{bib:arvesu-coussement-vanassche}.

Observe that despite formulas \eqref{eq:ch_mclassic} and
\eqref{eq:m1_mclassic}  corresponding to AT-systems of discrete
measures look similar to \eqref{eq:ch_rodrigues_util} and
\eqref{eq:m1_rodrigues_util} for Angelesco systems of measures,
there are significant differences between them: the backward
difference operator used in \eqref{eq:ch_mclassic} and
\eqref{eq:m1_mclassic} is not the same that the one involved in
\eqref{eq:ch_rodrigues_util} and \eqref{eq:m1_rodrigues_util}
(whenever $s(z)\not=z$, see formula \eqref{relation_delta_nabla}).
In addition, the mass points \eqref{eq:mass_points} in the
orthogonality relations \eqref{eq:ch_orth_cond_2} and
\eqref{eq:m1_orth_cond_2} differ from those used in the AT-systems
of discrete measures studied in
\cite{bib:arvesu-coussement-vanassche}. Note that the backward
difference operator \eqref{relation_delta_nabla} is defined base on
the mass points \eqref{eq:mass_points} on the $r$-star in order to
preserve the fundamental property of the Gamma function
$\Gamma(z+1)=z\Gamma(z)$, which plays an important role in our
construction. Another relevant fact, is the $\omega$-symmetry of the
system of measures (see Definition \ref{def:w-symmetric}). Indeed,
one can interpret the orthogonality conditions on the $r$-legged
star-like set \eqref{r-star} that define the Angelesco-polynomials
as rotated copies of the one given over the real interval.

Finally, the weak asymptotic of the two families of  Angelesco
polynomials studied here will be an interesting question. The
approach given in \cite{bib:aptekarev-arvesu} for discrete multiple
orthogonal polynomials should be explored for these families. In
addition, following the approach considered in this paper, other
Angelesco families of multiple orthogonal polynomials with respect
to discrete measures supported on an $r$-star should be
investigated.

\section*{Acknowledgements} Jorge Arves\'u would like to thank the Department of Mathematics at Baylor University for hosting his visit in Spring 2021 which stimulated this research. He also thanks Professor Andrei Mart\'{\i}nez-Finkelshtein, Baylor University, for his assistance and interesting discussions on Discrete Multiple Orthogonal Polynomials on star-like sets.

\noindent{\bf Funding.} J. Arves\'u acknowledges support by Agencia Estatal
de Investigaci\'on of Spain, grant number PID2021-122154NB-I00.



\begin{thebibliography}{9}

\bibitem{bib:abramowitz-stegun} M. Abramowitz and I. A. Stegun, {\it Handbook of Mathematical Functions},
National Bureau of Standards Applied Mathematics Series, 55, U.S. Government Printing Office, Washington, 1964.

\bibitem{bib:angelesco} A. Angelesco, Sur deux extensions des fractions continues alg\'ebriques,
Comptes Rendus Acad. Sci. Paris, 168, 262--265 (1919).

\bibitem{bib:aptekarev} A. I. Aptekarev, Multiple orthogonal polynomials, J. Comput. Appl. Math.
{\bf 99},  423--447 (1998).

\bibitem{bib:aptekarev-arvesu} A. I. Aptekarev and  J. Arves\'u, Asymptotics for multiple Meixner polynomials, J. Math. Anal. Appl. {\bf  411},  485--505 (2014).

\bibitem{bib:aptekarev-kalyagin-saff}  A. I. Aptekarev, V. A. Kalyagin, and  E. B. Saff, Higher-order three-term recurrences and asymptotics of multiple orthogonal
polynomials, Constr. Approx. {\bf 30}, 175--223 (2009).

\bibitem{bib:arvesu-coussement-vanassche} J. Arves\'u, J. Coussement,  and W. Van Assche, Some discrete multiple orthogonal polynomials, J. Comput. Appl. Math. {\bf 153},
19--45 (2003).

\bibitem{bib:arvesu-ramirez1} J. Arves\'u and  A. M. Ramirez-Abesturis, On the $q$-Charlier Multiple Orthogonal Polynomial, SIGMA {\bf 11},  026 (2015).

\bibitem{bib:arvesu-ramirez2} J. Arves\'u and  A. M. Ramirez-Abesturis, Multiple $q$-Kravchuk polynomials, Integral Transforms Spec. Funct. {\bf 32} (5-8),
361--376 (2020).

\bibitem{bib:arvesu-ramirez3} J. Arves\'u and  A. M. Ramirez-Abesturis, Multiple Meixner Polynomials on a Non-Uniform Lattice, Mathematics {\bf 8}(9),  1460 (2020).

\bibitem{bib:delvaux-lopez} S. Delvaux and  A. L\'opez, High order three-term recursions, Riemann--Hilbert minors and Nikishin systems on star-like sets, Constr. Approx. {\bf
37},  383--453 (2012).

\bibitem{bib:Gonchar} A. A. Gonchar, E. A. Rakhmanov, and  V. N. Sorokin, Hermite--Pad\'e approximants for systems of Markov-type functions, Mat. Sb. {\bf 188}, 33--58
(1997).

\bibitem{Hermite_exp} C. Hermite, {\it Sur la Fonction Exponentielle}, Gauthier-Villars, Paris, France (1874), 1--33.
\url{https://archive.org/details/surlafonctionexp00hermuoft/page/n1}

\bibitem{bib:kalyagin} V. A. Kalyagin, On a class of polynomials defined by two orthogonality relations, Mat. Sb. {\bf 110}, 609--627 (1979).

\bibitem{bib:kalyagin-ronveaux} V. A. Kalyagin and  A. Ronveaux, On a system of classical
polynomials of simultaneous orthogonality, J. Comput. Appl. Math.
{\bf 67}, 207--217 (1996).

\bibitem{bib:koekoek-lesky-swarttouw} R. Koekoek, P. A. Lesky, and  R. F. Swarttouw,
{\it Hypergeometric Orthogonal Polynomials and Their $q$-Analogs},
Springer-Verlag, Heidelberg-Berlin, 2010.

\bibitem{bib:lee} D. W. Lee, Classical multiple orthogonal polynomials of Angelesco
system, Appl. Numer. Math. {\bf 60},  1342--1351 (2010).

\bibitem{bib:leurs-vanassche2} M. Leurs and  W. Van Assche, Jacobi--Angelesco multiple orthogonal polynomials on an $r$-star, Constr. Approx. {\bf 51}, 353--381
(2020).

\bibitem{bib:leurs-vanassche1} M. Leurs and  W. Van Assche, Laguerre--Angelesco \
multiple orthogonal polynomials on an $r$-star, J. Approx. Theory
{\bf 250} 105324, 1--30 (2020).

\bibitem{bib:Mahler} K. Mahler, Perfect systems, Compos.
Math. {\bf 19}, 95--166 (1968). \url{http://eudml.org/doc/88959}

\bibitem{bib:miki-tsujimoto-vinet-zhedanov1}
H. Miki, S. Tsujimoto, L. Vinet, and  A. Zhedanov, An algebraic
model for the multiple Meixner polynomials of the first kind, J.
Phys.~A: Math. Theor. {\bf 45}, 325205 (2012).

\bibitem{bib:miki-tsujimoto-vinet-zhedanov2} H. Miki, L. Vinet, and  A. Zhedanov, Non-Hermitian oscillator Hamiltonians and multiple Charlier polynomials, Phys.
Lett. A {\bf 376}, 65--69 (2011).

\bibitem{bib:ndayiragije1-vanassche} F. Ndayiragije1 and  W. Van Assche, Multiple Meixner
polynomials and non-Hermitian oscillator Hamiltonians, J. Phys. A:
Math. Theor. {\bf 46}, 505201, 17 p. (2013).

\bibitem{nu} A. F. Nikiforov and  V. B. Uvarov, {\it Special Functions of Mathematical
Physics. A Unified Introduction with Applications}, Birkh\"{a}user,
Basel, 1988.

\bibitem{bib:nikiforov-uvarov-suslov} A. F. Nikiforov, S. K. Suslov, and  V. B. Uvarov, {\it Classical Orthogonal
Polynomials of a Discrete Variable}, Springer-Verlag, Berlin, 1991.

\bibitem{bib:nikishin} E. M. Nikishin, A system of Markov functions, Moscow Univ. Math.
Bull. {\bf 34}, 63--66 (1979).

\bibitem{bib:kn_Nikishin} E. M. Nikishin, On simultaneous Pad\'{e} approximations, Mat.
Sb. {\bf  113}, 499--519 (1980) (in Russian); English transl. Math.
USSR Sb. {\bf 41} (1982).

\bibitem{bib:nikishin-sorokin} E. M. Nikishin and  V. N. Sorokin, {\it Rational Approximation and Orthogonality}, Translations of Mathematical Monographs, 92,
American Mathematical Society, Providence, RI, 1991.

\bibitem{bib:Sorokin3}
V. N. Sorokin, Simultaneous Pad\'e approximants for finite and
infinite intervals, Izv. Vyssh. Uchebn. Zaved. Mat. (8), 45--52
(1984). \url{http://mi.mathnet.ru/eng/ivm8510}

\bibitem{bib:sorokin2} V. N. Sorokin, A generalization of Laguerre polynomials and
convergence of simultaneous Pad\'e approximants, Uspekhi Mat. Nauk
{\bf 41}, 245--246 (1986).

\bibitem{bib:sorokin1} V. N. Sorokin, A generalization  of classic orthogonal
polynomials and the convergence of simultaneous Pad\'e
approximants, J. Soviet Math. {\bf 45} (6), 1461--1499 (1989).

\bibitem{bib:vanassche} W. Van Assche, Nearest neighbor recurrence relations
for multiple orthogonal polynomials, J. Approx. Theory {\bf 136},
1427--1448 (2011).

\bibitem{bib:vanassche-coussement} W. Van Assche and  E. Coussement, Some classical
multiple orthogonal polynomials, J. Comput. Appl. Math. {\bf 127},
317--347 (2001).

\end{thebibliography}
\end{document}